\documentclass[draft,11pt,leqno]{amsart}

\usepackage{amsmath,amssymb,amscd,amsthm}
\usepackage{latexsym}
\usepackage[dvips]{graphics,epsfig}

\newtheorem{theorem}{Theorem}
\newtheorem{lemma}{Lemma}
\newtheorem{proposition}{Proposition}

\newtheorem{definition}{Definition}
\newtheorem{corollary}{Corollary}

\newtheorem{claim}{Claim}

\newcommand{\f}[2]{\frac{#1}{#2}}
\newcommand{\dpr}[2]{\left\langle #1,#2 \right\rangle}

\newcommand{\al}{\alpha}

\newcommand{\de}{\delta}

\newcommand{\De}{\Delta}
\newcommand{\ve}{\varepsilon}

\newcommand{\ka}{\kappa}
\newcommand{\la}{\lambda}

\newcommand{\si}{\sigma}

\newcommand{\vp}{\varphi}
\newcommand{\om}{\omega}

\newcommand{\rd}{{\mathbf R}^d}

\newcommand{\rone}{\mathbf R^1}

\newcommand{\cl}{\mathcal L}

\newcommand{\cg}{\mathcal G}

\newcommand{\ca}{\mathcal A}

\newcommand{\cm}{\mathcal M}
\newcommand{\ct}{\mathcal T}

\newcommand{\cb}{\mathcal B}

\newcommand{\cd}{\mathcal D}

\newcommand{\cz}{\mathcal Z}

\newcommand{\cc}{\mathcal C}
\newcommand{\ch}{\mathcal H}

\newcommand{\p}{\partial}

\newcommand{\beq}{\begin{equation}}
\newcommand{\eeq}{\end{equation}}
\newcommand{\beqna}{\begin{eqnarray*}}
\newcommand{\eeqna}{\end{eqnarray*}}
\newcommand{\beqn}{\begin{equation*}}
\newcommand{\eeqn}{\end{equation*}}
\newcommand{\bp}{\begin{proof}}
\newcommand{\ep}{\end{proof}}
\newcommand{\bprop}{\begin{proposition}}
\newcommand{\eprop}{\end{proposition}}
\newcommand{\bt}{\begin{theorem}}
\newcommand{\et}{\end{theorem}}
\newcommand{\bex}{\begin{Example}}
\newcommand{\eex}{\end{Example}}
\newcommand{\bc}{\begin{corollary}}
\newcommand{\ec}{\end{corollary}}
\newcommand{\bcl}{\begin{claim}}
\newcommand{\ecl}{\end{claim}}
\newcommand{\bl}{\begin{lemma}}
\newcommand{\el}{\end{lemma}}

\begin{document}

\title
[Stability and instability of traveling waves of second order in time PDE's]
{Linear stability analysis for traveling waves of second order in time PDE's}

\author{Milena Stanislavova} 

\author{Atanas Stefanov}

\address{Milena Stanislavova
Department of Mathematics, 
University of Kansas, 
1460 Jayhawk Boulevard,  Lawrence KS 66045--7523}
\email{stanis@math.ku.edu}
\address{Atanas Stefanov 
Department of Mathematics, 
University of Kansas, 
1460 Jayhawk Boulevard,  Lawrence KS 66045--7523}
\email{stefanov@math.ku.edu}

\thanks{
Stanislavova  supported in part by  
NSF-DMS 0807894. 
Stefanov supported in part by  
NSF-DMS \# 0908802 .}
\date{\today}



\begin{abstract}
We study traveling waves $\vp_c$ of second order in time PDE's  $u_{tt}+\cl u+N(u)=0$. The linear stability analysis for these models is reduced to the question for stability of quadratic pencils in the form  $\la^2Id+2c\la \p_x+\ch_c$, where $\ch_c=c^2 \p_{xx}+\cl+N'(\vp_c)$. 

If $\ch_c$ is a self-adjoint operator, with a simple negative eigenvalue and a simple eigenvalue at zero, then   we completely characterize the linear stability of $\vp_c$. More precisely, we introduce an explicitly computable index $\om^*(\ch_c)\in (0, \infty]$, so that the wave $\vp_c$ is stable if and only if $|c|\geq \om^*(\ch_c)$. 
The results are applicable both in the periodic case and in the whole line case. 

The method of proof involves a delicate analysis of a  function $\cg$, associated with $\ch$, whose positive zeros are exactly the positive (unstable) 
eigenvalues of the pencil $\la^2Id+2c\la \p_x+\ch$.  We would like to emphasize that the function $\cg$ is not the  Evans function for the problem, but rather a new object that we define herein, which fits the situation rather well. 

As an application, we consider three classical 
  models - the ``good'' Boussinesq equation,  
 the Klein-Gordon-Zakharov (KGZ) system and the fourth order  beam equation. In the whole line case, for the Boussinesq case and the KGZ system (and as a direct application of the main results), we compute explicitly the set of speeds which give rise to linearly stable traveling waves (and for all powers of $p$ in the case of Boussinesq). This result is new for the KGZ system, while it generalizes the results of \cite{AS} and \cite{Alex},  which apply to the case $p=2$. 
 For the beam equation, we provide an explicit formula 
 (depending of the function $\|\vp_c'\|_{L^2}$), which works for all $p$ and for both the periodic and the whole line cases. 
 
 Our results complement (and exactly match, whenever they exist) the results of a long line of investigation regarding the related notion of orbital stability of the same  waves. 
Informally, we have found that in all the examples that we have looked at, 
our theory goes through, whenever the Grillakis-Shatah-Strauss theory applies. We believe that the  results in this paper (or a variation thereof) will enable the linear stability analysis as well as asymptotic stability analysis for  most models in the form $u_{tt}+\cl u+N(u)=0$. 
 
\end{abstract}

\maketitle

\section{Introduction and motivation}
The main motivation of our study is the following abstract second order in time  nonlinear PDE
\begin{equation}
\label{abst}
u_{tt}+ \cl_x u+N(u)=0, \ \ (t,x)\in \rone_+\times \rd\ \  \textup{or}\ \ 
(t,x)\in \rone\times [-L,L]^d, 
\end{equation}
where $\cl_x$ is a given linear operator, acting on the $x$ variable and $N(u)$ is the nonlinear term. We outline some relevant examples of interest, see Section \ref{sec:examples} below, but notice that we do consider both periodic boundary conditions as well as vanishing at infinity solutions. These two scenarios will be considered simultaneously, since our method works equally well in both cases. 

Our interest is in the study of the stability properties of traveling waves, that is, solutions in the form 
$\vp(x+\vec{c} t)$. Clearly, these  satisfy the stationary PDE 
\begin{equation}
\label{abs1}
\cl_x \vp+\sum_{i,j=1}^d c_i c_j \vp_{x_i}\vp_{x_j}+N(\vp)=0
\end{equation}
We take the ansatz $u=\vp(x+\vec{c} t)+v(t, x+\vec{c} t)$ and plug it in \eqref{abst}. Taking into account 
\eqref{abs1} and dropping all quadratic and higher order terms in $v$, we arrive at the following 
\begin{equation}
\label{abs2}
v_{tt}+2 \dpr{\vec{c}}{\nabla_x v_t}+\sum_{j,k=1}^d c_j c_k \f{\p^2 v}{\p x_j \p x_k}+\cl v+
N'(\vp)v=0
\end{equation}
Therefore, introducing the operator $H_{\vec{c}}=\cl_x+ 
\sum_{j,k=1}^d c_j c_k  \p_{x_j} \p_{x_k}+ N'(\vp)$, we are lead to study  the following problem 
\begin{equation}
\label{abs3}
v_{tt}+2 \dpr{\vec{c}}{\nabla_x v_t}+H u=0.
\end{equation}
This is the linearized stability problem for the nonlinear equation \eqref{abst} in a vicinity of the special solution $\vp(x+\vec{c} t)$. Before we move on with our exposition, we provide some informal definitions and discussions to motivate our interest in the problem. There are several notions of stability that are of interest. Generally, linear stability is easier to check than (and is often a prerequisite to) nonlinear stability. Even within the linear stability, we distinguish between spectral and linear stability. Namely, for an evolutionary problem in the form $z_t=\cm z$, where $\cm$ is a closed operator generating $C_0$ semigroup, we say that we have {\it spectral stability}, if $\si(\cm)\subset \cz_-=\{\la:  \Re\la\leq 0\}$.  We say that the same problem has {\it linear stability}, if the solutions grow at infinity slower than any exponential\footnote{Here, note that one cannot require bounded orbits. Indeed,  even in the ODE case, if one takes a Jordan block $A$ of size $d$,  $e^{t A}$ has polynomial growth $t^d$}, i.e. for every $\de>0$, 
$\lim_{t\to \infty} e^{- \de t}\|z(t)\|=0$. The relationship between spectral and linear stability has been well-explored and documented in the literature, 
and we will not dwell on it, except to point out that in principle (and in the presence of the so-called spectral mapping theorem for the generator $\cm$), these are equivalent and amount to lack of exponentially growing modes, that is solutions in the form $\cm \psi=\la \psi$ for $\la: \Re\la>0$. 

On the other hand, we have two distinct notions of nonlinear stability - orbital and asymptotic stability. Assuming for simplicity that the only invariance of the system is translation, orbital stability requires that a solution for the full nonlinear equation that starts close to a   traveling wave stay close for all times 
to a (time-dependent) translate of  the starting wave. Asymptotic stability requires a bit more, namely that the perturbed profile will actually converge to a (time-dependent) translate of the wave. 

There is a large body of literature that deals with this problem in 
various models. We would like to point out that the powerful methods of Grillakis-Shatah-Strauss, \cite{GSS1}, \cite{GSS2} reduce (in most cases)
 the problem of orbital stability to checking certain conditions on the linearized functionals. We mention in this regard the papers \cite{Bona}, \cite{Chen}, \cite{Liu1}, \cite{Lev}, which treat models of interest to us in this paper. 
 We also   note that establishing  orbital instability for a given problem seems to be harder and requires more problem specific efforts, \cite{Liu1}. 
 
 Our interest here is the question of linear stability of such traveling waves 
$\vp_c$, that is, whether there are exponentially growing solutions of \eqref{abs3} in the form $e^{\la t} \psi(x)$.  We pose the following

 {\bf Question:} 
For which values of $\vec{c}$, the corresponding traveling wave $\vp_c$ determined by \eqref{abs1} is linearly/spectrally stable? More precisely, for which $\vec{c}$, the equation \eqref{abs3} has a solution in the form $e^{\la t}\psi$? \\

Note that this question has been studied  thoroughly in the last twenty years. Here, the methods are completely different, although the results in the end must be related with those obtained by orbital stability methods. After all one expects orbitaly stable waves to be linearly stable and vice versa (although this sometimes fails at points where the stability character changes). 

The Evans function method has been used in connection to the linear stability of such waves, \cite{AS}, \cite{PW1}, \cite{PW},  \cite{LL}. However, the method has well-known limitations, in particular dealing with  systems, which is the issue at hand here, since the equations are second order in time. In fact, the   papers quoted here are among the few that deal with second order in time equations. 

Another method that has been developed is the method of ``indices counting'', which basically relates the number of unstable eigenvalues of self-adjoint entries of $\cm$ to the number of unstable eigenvalues of $\cm$ itself. This has been mostly useful in the KdV and Schr\"odinger type systems, but it has also played important role in general Hamiltonian and dissipative systems,  \cite{KKS1}, \cite{KKS2}, \cite{Kap1}, \cite{CP},  \cite{DK}. We would like to point out that some of these results have enabled the consideration of spatially periodic waves, \cite{DK}, which is one of the goals of this project as well. 

A third method for establishing mostly sufficient conditions for instability is through a direct construction of unstable modes, \cite{Liu1}, 
\cite{Lin}. In fact, there are many other impressive results in the literature, mostly for standing waves, which provide ``instability by blow up'' for such waves. These are done typically by constructing  clever Lyapunov functionals, associated with the unstable modes.  

Our goal in this paper is to develop a fairly general and systematic theory, which treats the question for linear stability (i.e. the existence/nonexistence of such $\psi$)  of  traveling waves of second order in time systems. {\it  Indeed, we present a complete answer to this question, in the case of \eqref{abs3}, when $H$ is self-adjoint, with at most one negative eigenvalue and $d=1$.  }

\subsection{Examples}
\label{sec:examples}
We consider three basic examples that fit this category, although there are numerous others to which our results will be applicable. The reader should also bear in mind that one might consider both periodic and whole line boundary conditions, for the 
models that we list below. To keep the discussion simple, we will mostly stick to the whole line case, the periodic cases will be addressed in a subsequent publication, \cite{HSS}.  

\subsubsection{The ``good'' Boussinesq models} Our first example   is   the ``good' Boussinesq-type  model, 
\begin{equation}
\label{b:2}
u_{tt}+u_{xxxx}-u_{xx}+(u^p)_{xx}=0, \ \ (t,x)\in \rone_+\times \rone 
\end{equation} 
This model was considered by Bona-Sachs, \cite{Bona} as a model for propagation of small amplitude, long waves on the surface of water. This is indeed an equation that belongs to a family of Boussinesq  models, which all have the same level of formal validity, however \eqref{b:2} exhibits some desirable features, like local   well-posedness in various function spaces, \cite{Bona}. Interestingly, global well-posedness for \eqref{b:1} does not hold, even if one requires smooth initial data with compact support. In fact, there are  ``instability by blow-up'' results for such unstable traveling waves for this equation. 

It is easy to see that there exists one-parameter family of traveling waves of the form $\vp(x-c t), |c|\in (-1,1)$, which obey the equation 
\begin{equation}
\label{s:1}
c^2\vp+\vp''-\vp+\vp^p=0
\end{equation}
and which have the explicit form 
$$
\vp_c(\xi)=\left[\left(\f{p+1}{2}\right) (1-c^2)  \right]^{\f{1}{p-1}} 
 sech^{\f{2}{p-1}}\left(\f{\sqrt{1-c^2} (p-1)}{2} \xi\right).
$$
One should of course recognize the standard $sech$ solitons of the generalized KdV hierarchy. This is not surprising, since the governing ODE for the traveling waves of  \eqref{b:1} is related to the corresponding ODE for gKdV. 

In \cite{Bona}, the authors  proved orbital stability for these solutions, provided   $p<5, \f{\sqrt{p-1}}{2}<|c|<1$, while nonlinear instability was established in \cite{Liu1}. Both of these results were  achieved via the methods of the Grillakis-Shatah-Strauss theory (developed in \cite{GSS1}, \cite{GSS2})   for orbital stability/instability. 
\subsubsection{The Klein-Gordon-Zakharov system} 
Consider the following system of coupled wave equations\footnote{Here the coefficient $\f{1}{2}$ in front of the $ (|u|^2)_{xx}$ can be taken to be arbitrary by rescaling, but the particular choice of $\f{1}{2}$ will be convenient in the spectral analysis}
\begin{equation}
\label{KGZ}
\left|
\begin{array}{l}
u_{tt}-u_{xx}+u + n u=0\ \  (t,x)\in \rone_+\times \rone \\
n_{tt}-n_{xx}-\f{1}{2}(|u|^2)_{xx}=0,
\end{array}
\right.
\end{equation}
which describes the interaction of a Langmuir wave and an ion acoustic wave in a plasma\footnote{namely, the complex $u$ describes the fast scale component of the electric field, while the real valued $n$ stands for the deviation of ion density.}, \cite{Dendy}. Taking the traveling wave ansatz $u(t,x)=\vp(x-ct), n(t,x)=\psi(x-ct)$, we derive from the second equation that $(c^2-1) \psi''=\f{1}{2}(\vp^2)''$, whence 
$(c^2-1)\psi= \f{\vp^2}{2}$, since the functions $\vp, \psi$ decay at infinity. Using this relation in the first equation yields the following ODE for $\vp$
\begin{equation}
\label{eq:k1}
-(1-c^2)\vp''+\vp-\f{\vp^3}{2(1-c^2)}=0
\end{equation}
It is well-known (and in fact can be checked easily, based on simple rescaling arguments in the ODE's governing traveling wave solutions for the KdV equation) that \eqref{KGZ} admits an one parameter family of traveling wave solutions in the form $u(t,x)=\vp(x-ct), n(t,x)=\psi(x-ct)$ for $c\in (-1,1)$, where 
\begin{equation}
\label{sol:KGZ}
\left|
\begin{array}{l}
  \vp(y)=  2 \sqrt{1-c^2} sech\left(\f{y}{\sqrt{1-c^2}}\right) \\ 
 \psi(y)=-2 sech^2 \left(\f{y}{\sqrt{1-c^2}}\right).
 \end{array} 
\right.
\end{equation}
In fact, there exists a two-parameter system  of solitary-traveling waves of which the family displayed above is a particular case, \cite{Chen}. In the same paper,  the author shows orbital stability for such waves, provided $1>|c|>\f{\sqrt{2}}{2}$. Below, we solve completely the question for linear stability.  Namely, for $|c|\in [\f{\sqrt{2}}{2},1)$ we have linear stability, whereas for $|c|\in [0, \f{\sqrt{2}}{2})$ we have linear instability.

\subsubsection{The nonlinear beam equation} 
Another relevant example to consider is the so-called  beam equation 
\begin{equation}  
\label{b:1}
u_{tt}+\De^2 u+u-|u|^{p-1} u=0, \ \ (t,x)\in \rone\times \rd\ \textup{or}\ \ 
(t,x)\in \rone\times [-L,L]^d, 
\end{equation}
where $p>1, L>0$ and we either require periodic boundary conditions (in the case $x\in [-L,L]$) or 
vanishing at infinity for $x\in \rd$. 

This  equation  has been studied extensively in the literature. It seems that the earliest work on the subject goes back to \cite{MW}, where \eqref{b:1} is proposed as a model of a suspension bridge. In this paper, the authors have also showed  the existence of some traveling wave solutions, that is, solutions in the form $\vp(x+\vec{c} t)$. 

In \cite{Lev} it was proved that traveling wave solutions exist, in the whole space context, whenever\footnote{as we shall see in our arguments later on, this is a natural space of parameters} $|\vec{c}|\in (0, \sqrt{2})$. 
The proof is variational in nature, and this allowed the author to derive conditions for orbital stability/instability of such  waves for different speeds. 
In short, the conclusion was that such traveling waves are orbitally 
unstable for small speeds, while an orbitally stable solutions are observed for values of the parameter $|\vec{c}|\sim \sqrt{2}, |\vec{c}|<\sqrt{2}$. 
Various other works have explored the decay and scattering 
 properties of the linear beam equation 
$u_{tt}+\De^2+u=0$, \cite{Lev2, Miao,Pau}.  

\subsection{Some general features of the operators $H$} 
In this section, we construct the operators $H$ that arise in \eqref{abs3} and discuss some of their general features. Of particular importance will be the fact that they are self-adjoint and the fact that they have exactly one negative eigenvalue. 

\subsubsection{The beam equation - spectral picture} 
\label{sec:1.2.1}
 Clearly, such a function $\vp$ will satisfy the ODE 
\begin{equation}
\label{b:5}
\De^2 \vp +\sum_{i,j=1}^d c_i c_j  \vp_{x_i x_j}+\vp-|\vp|^{p-1}\vp=0, 
\ \ x\in (-L,L)^d,
\end{equation}
which is supplied by the usual periodic boundary conditions at $-L, L$ (or vanishing at infinity for the whole line case). 
Note that we do not necessarily look for positive solutions of \eqref{b:5},  
nor do we expect to find any in general (because of the lack of maximum principle for $\De^2$, the ground state are not positive), thus the presence of the absolute value in the non-linear term.  In the one dimensional case, which will be the main subject of our investigation, \eqref{b:5} reduces to 
\begin{equation}
\label{b:6}
c^2 \vp''+\vp'''' +\vp-|\vp|^{p-1}\vp=0, 
\ \ -L<x<L. 
\end{equation}
Clearly, the question for the existence of such $2L$ periodic solutions $\vp$ is a 
non-trivial one, but we refer the reader to the upcoming publication \cite{DS}, where it will be addressed. Note though that if $\vp$ is constructed in the space $L^2_{per.}[-L,L]$, then one can  automatically bootstrap its smoothness to infinity, thanks to the form of \eqref{b:6}. 

Even less obvious issue is the  linear stability of such waves. In order to simplify matters, let us assume henceforth that $d=1$ and $p$ is an odd integer (in particular, one should keep in mind the model case $p=3$). 
Now, for any such periodic wave $\vp$ satisfying \eqref{b:6}, set 
$u(t,x)=\vp(x+c t)+v(t,x+ct)$, which we plug in \eqref{b:1}. 
After ignoring all the quadratic terms in $v$, taking into account \eqref{b:6} and changing variables $x- ct\to x$, we obtain the following {\it linearized equation} around the traveling wave profile $\vp$ 
 \begin{equation}
\label{b:10}
v_{tt}+2 c v_{t x}+ v_{xxxx}+c^2 v_{xx}+v- p \vp^{p-1} v=0, 
\ \ -L<x<L. 
\end{equation}
This is then a proper time to introduce the self-adjoint operator $\ch$ 
\begin{equation}
\label{b:12}
\ch:= \p_x^4+c^2 \p_x^2+1- p \vp^{p-1},
\end{equation} 
with the domain $D(\ch):=H^4_{per.}[-L,L]$. 
Note that $\ch=\ch(c,p,\vp)$ depends implicitly on the speed $c$, the parameter $p$ and  the wave $\vp$ itself. 

Another fact that should be stated at this point is that $\ch$ has zero on its spectrum, with eigenfunction $\vp'$. Indeed, taking a spatial derivative in \eqref{b:6} (under the assumption that $p$ is an odd integer) 
implies $\ch[\vp']=0$. Finally,  the bottom of the spectrum of $\ch$ is a negative number\footnote{note that Weyl's theorem, in the whole line case 
 $\si_{a.c.}(\ch)=[1, \infty)$}. Indeed, 
$$
\dpr{\ch \vp}{\vp}=-(p-1)\int_{-L}^L \vp^{p+1}(x) dx<0,
$$
hence $\si(\ch)\cap (-\infty, 0)\neq \emptyset.$  
Unfortunately, we cannot verify that this eigenvalue is simple, although this is easily seen in numerical simulations, for all values of the parameter $c$. 

Note that the computations above for the beam equation in the periodic case apply equally well for the whole line homoclinic solution produced by Levandosky, \cite{Lev}. 

\subsubsection{The ``good'' Boussinesq model: spectral picture}
Here we outline the relevant operators for the Boussinesq model.   As we shall see though, the scheme outlined in \eqref{abs3} does not quite work and one needs a little trick to make it work. 

In any case,  if we impose the ansatz $u=\vp_c(x+c t)+v(t, x+ct)$ and we get the equation $v_{tt}+2 c v_{t x}+ T v=0$, where 
$$
T v=\p_x^4 v -(1-c^2)\p_x^2 v +p (\vp_c^{p-1} v)_{xx}
$$
Note that the operator $T$ is not self-adjoint and thus unsuitable for our method. However, if we introduce the variable $z: z_x=v$, we get the following linearized equation in terms of $z$, 
\begin{equation}
\label{zz:2}
z_{ttx}+2 c z_{t xx}+ T[z_x]=0.
\end{equation}
Note that $T[z_x]=\p_x[H[z]]$, where 
\begin{equation}
\label{z:2}
H z=\p_x^4 z-(1-c^2)\p_x^2 z +p (\vp_c^{p-1} z_x)_{x}
\end{equation}
Thus, the linearized equation becomes $\p_x[z_{tt}+2 c z_{tx}+ H z]=0$, which in view of the conditions $\lim_{|x|\to \infty} |\p_x^\al z(x)|=0$, implies $z_{tt}+2 c z_{tx}+ H z=0$. Thus, the linearized equation to consider is again $z_{tt}+2 c z_{tx}+H z=0$, where $H$ is defined in \eqref{z:2}. As we prove later on, $H$   has  one simple eigenvalue at zero and one simple negative eigenvalue. 

{\it Note that while it is clear the presence of unstable mode  for $z_{tt}+2 c z_{tx}+ H z=0$ implies linear instability for \eqref{zz:2} (via $v=z_x$), the reverse it is not immediately clear.}  The reason for that is that 
the operator $\p_x$ is not invertible and one may potentially have a solution 
$e^{\la t} V(x)$ of \eqref{zz:2}, while the corresponding ``solution'' 
 $e^{\la t}Z(x)=e^{\la t}\p_x^{-1} V(x)$ may not be an $L^2$ function.

\section{Main Results}
In this section, we present our results. We start first with a abstract form of the linearized problem, which involves the theory of quadratic pencils.  

\subsection{Stability/instability results for quadratic pencils}
In connection with \eqref{b:10}, we may write it schematically in the form
$$
v_{tt}+2 c v_{t x}+\ch v=0.
$$
\begin{definition}
\label{defi:1}
We say that the periodic wave $\vp$ is linearly unstable, if there exists  
$\la: \Re \la>0$, and a function $\psi$, so that 
then the following equation is satisfied
\begin{equation}
\label{b:15}
\la^2 \psi+2c\la \psi_x+ \ch\psi=0
\end{equation}
Otherwise, we say that the traveling wave $\vp$ is linearly stable. 
\end{definition}
In order to study the stability of such waves, we actually consider the much 
 more general framework  of stability/instability of operator pencils of the form 
$$
L(\la)=\la^2 Id +2 c \la \p_x+ H,
$$
where $H$ is a fixed self-adjoint operator, acting on $L^2$, which satisfies 
certain mild constraints. This question has been considered before, 
mainly in a pure operator-theoretic sense. Note the work   Chugunova-Pelinovsky, \cite{CP1} in which they consider similar problems, but where the operator $H$ is a finite symmetric matrix. In our work, we will need the general theory developed in the works of Shkalikov, \cite{Shkal}, see also Azizov-Iokhvidov, \cite{AI} for a thorough review of the available results.   We also 
develop some of the theory below, as we shall need it for our purposes, see Section \ref{sec:3.1}  below.  

Going back to the problem at hand,  we consider linear, second-order in time equations in the general form 
\begin{equation}
\label{g1}
u_{tt}+2 \om u_{tx}+H u=0, (t,x)\in \rone\times \rone \ \ \textup{or} \ \ \rone\times[-L,L]
\end{equation}
where $H=H_c$ is a self-adjoint operator acting on $L^2$, 
with domain $D(H)$ and $\om$ is a real parameter.  Note that it is better at this point to consider $\om$ as an independent parameter, and to ignore the fact that in the applications $\om=c$. 

 Based on our examples, we saw that it is reasonable  to assume that the self-adjoint operator $H$ has one simple negative eigenvalue, a  simple eigenvalue at zero (which is naturally generated by $\vp'$) and the rest of the spectrum is contained in $(0,\infty)$, with a spectral gap. Even though our approach is quite systematic and should be able to handle a more general situation\footnote{for example multiple eigenvalues at zero, which is the case in higher dimension}, we choose to implement these assumptions for simplicity of the exposition. Thus, we require 
\begin{equation}
\label{A} 
\left\{ \begin{array}{l} 
\si(H)=\{-\de^2\} \cup \{0\}\cup \si_+(H), \si_+(H)\subset (\si^2, \infty), \si>0\\
 H \phi=-\de^2 \phi,  \ \dim[Ker(H+\de^2)]=1 \\
  H \psi_0=0, \ \  \dim[Ker(H)]=1\\
  \|\phi\|=\|\psi_0\|=1
\end{array}
\right.
\end{equation}
 We shall  need the following spectral projection operators 
\begin{eqnarray*}
& & P_0: L^2\to \{\phi\}^\perp; P_0 h:=h-\dpr{h}{\phi}\phi \\
& & P_1: L^2\to \{\phi, \psi_0\}^\perp; P_1 h:=h-\dpr{h}{\phi}\phi- 
\dpr{h}{\psi_0}\psi_0
\end{eqnarray*}
Our next assumption is essentially that $H_1$ is of order higher than  one. We put it in the following  form: for all $\tau>>1$, we require 
\begin{equation}
 \label{E}
   (H+\tau)^{-1/2} \p_x (H+\tau)^{-1/2}, (H+\tau)^{-1} \p_x \in \cb(L^2),  
 \end{equation}
 Note that the quantities in \eqref{E} are well-defined, since for all $\tau>>1$, we have  that $H+\tau\geq (\tau-\de^2)Id>0$. 
 
 An easy consequence of \eqref{E} is that 
 $H^{-1} P_1\p_x P_1\in \cb(L^2)$ as well. This follows easily from the resolvent identity, since $H^{-1} P_1\p_x P_1=P_1 (H+\tau)^{-1} \p_x P_1+\tau H^{-1} P_1 (H+\tau)^{-1} \p_x P_1$. In addition, the  following non-degeneracy condition  is also required 
 \begin{equation}
\label{B}    
\dpr{\phi'}{\psi_0}\neq 0.
\end{equation}
Lastly, we assume that $H$ has real coefficients. We formulate as follows 
\begin{equation}
\label{G}
\overline{H h}=H\bar{h}.
\end{equation} 

An important observation, that we would like to make right away (and which will be used repeatedly in our arguments later on) is that for every $\la>0$, the operator $(H+\la^2):\{\phi\}^\perp\to \{\phi\}^\perp$ is invertible, since 
$(H+\la^2)|_{\{\phi\}^\perp}\geq \la^2 Id$.

The following theorem is the main result of this paper. 
\newpage  
\begin{theorem} 
\label{theo:5}
 Let $H$  be a self-adjoint operator on $L^2$. Assume that  it satisfies the structural assumption \eqref{A}, \eqref{E} as well as the 
 non-degeneracy assumption  \eqref{B} and  \eqref{G}.  
 
Then, if $\dpr{H^{-1}[\psi_0']}{\psi_0'}\geq  0$, we have 
instability (in the sense of Definition \ref{defi:1}) for all values of $\om\in \rone$.

 Otherwise, supposing $\dpr{H^{-1}[\psi_0']}{\psi_0'}<0$, 
we have 
 \begin{itemize}
 \item 
  the problem \eqref{g1} is unstable if 
   $\om$ satisfies the inequality
 \begin{equation}
 \label{m:2}
 0\leq |\om|< 
 \frac{1}{2\sqrt{-\dpr{H^{-1}[\psi_0']}{\psi_0'}}}=:\om^*(H)
 \end{equation}
 \item  the problem \eqref{g1} is stable, if 
   $\om$ satisfies the reverse inequality 
   \begin{equation}
 \label{m:3}
  |\om|\geq  \om^*(H)
 \end{equation}
 \end{itemize} 
  
\end{theorem}
 {\bf Remarks: } 
 \begin{itemize}
 \item Note that the result is  a complete characterization of the stability and instability properties of the abstract quadratic pencil problem \eqref{g1}. In essence, it says that there is a critical number $\om^*(H)$,  below which the problem is unstable and above which, there is stability. 
  \item Note that the critical number $\om^*(H)$ may also depend on\footnote{and in fact, for the applications in mind, $H$ will have the form of \eqref{b:12}, where the dependence on $c$  is pretty transparent} $c$, so that the equality   $c=\om^*(H,c)$ (as in \eqref{b:15}) will be an implicit one. In any case, the conclusions \eqref{m:2} and \eqref{m:3} hold true, and should be used to determine stability,  regardless of this implicit dependence. 
  \item A more succinct way of defining $\om^*(H)$ is the following 
  $$
  \om^*(H)=\left\{\begin{array}{cc}
  +\infty & \textup{if}\ \ \dpr{H \psi_0'}{\psi_0'}\geq 0 \\
  \frac{1}{2\sqrt{-\dpr{H^{-1}[\psi_0']}{\psi_0'}}} & \textup{if}\ \  \dpr{H \psi_0'}{\psi_0'}<0 
  \end{array}
  \right.
  $$ 
 Then, we can characterize stability in terms of the inequality as follows: $|\om|\geq \om^*(H)$. 
  \end{itemize} 
  An immediate corollary is the following ``spectral stability/instability'' type of statement. We can write \eqref{g1} in the form 
  $$
  \left(\begin{array}{c} u \\ u_t   \end{array}   \right)_t=
  \left(\begin{array}{cc} 0 &  1  \\
  -H & -2\om\p_x  \end{array}   \right) \left(\begin{array}{c} u \\ u_t   \end{array}   \right)=:\ct \left(\begin{array}{c} u \\ u_t   \end{array}   \right)
  $$
  \begin{corollary}
  \label{cor:l1}
  In the statement of Theorem \ref{theo:5}, assume in addition that 
  $[H h(-\cdot)](x)=(H h)(-x)$. Then, in the cases of instability, there is $\la>0$, so that $\la, -\la$ are both eigenvalues of $\ct$ and moreover 
  $$
  \si(\ct)\subset \{\la\}\cup \{-\la\} \cup i\rone. 
  $$
  If on the other hand, there is stability, we have $\si(\ct)\subset   i\rone. $ 
  \end{corollary}
   
  \subsection{Applications}
  
  \subsubsection{Stability for traveling waves of the ``good'' Boussinesq equation}
  The linear stability/instability,  for the case $p=2$, was considered in the literature. The linear instability for  $|c|<\f{\sqrt{p-1}}{2}$ was shown by Fal'kovich-Spector-Turitsyn,  \cite{FST}. A more complete analysis, which handled both the stability and the instability regimes was performed by Alexander and Sachs, \cite{AS}. The idea behind the proof is based on an  Evans function calculations, but the authors needed some computer assistance\footnote{We were also informed, \cite{Alex}, that Alexander, Sachs and Pego had a different proof, which is unpublished. }. 
  \begin{theorem}
  \label{theo:bous}
  The traveling wave $\vp_c$ of the Boussinesq equation \eqref{b:1} is linearly unstable, if $p\geq 5$. If $2\leq p<5$, then it is linearly unstable if 
  $0\leq |c| < \f{\sqrt{p-1}}{2}$ and linearly stable, when $\f{\sqrt{p-1}}{2}\leq |c|<1$.  
  \end{theorem}
  {\bf Remarks:} Note that the result in Theorem \ref{theo:bous} matches precisely the results for orbital stability for the ``good'' Boussinesq model. Recall that orbital stability was shown by Bona-Sachs for $|c|>\f{\sqrt{p-1}}{2}$ and the nonlinear instability by Liu, \cite{Liu1}.  
  
  The critical case $|c|=\f{\sqrt{p-1}}{2}$ provides a (marginally) linearly stable situation, according to our result, but note that Liu, \cite{Liu1} shows nonlinear instability there. This is due to the fact that in this particular case, there is an additional eigenvalue at zero (which is unaccounted for in terms of symmetries), which is responsible for a secular nonlinear instability. 
  \subsubsection{Stability for the Klein-Gordon-Zakharov system} 
  \begin{theorem}
  \label{theo:KGZ} 
 Let $c\in (-1,1)$. Then, the traveling wave solution $(\vp(x-ct), \psi(x-ct))$ described in \eqref{sol:KGZ} is spectrally/linearly stable for 
 $|c|\in [\f{\sqrt{2}}{2},1)$ and linearly/spectrally unstable for 
 $|c|\in [0, \f{\sqrt{2}}{2})$. 
  \end{theorem}
  {\bf Remarks:} Note that the linear stability results match precisely the orbital stability results in \cite{Chen}, except at the endpoints $|c|=\f{\sqrt{2}}{2}$. At this point, we have linear stability, according to Theorem \ref{theo:KGZ}, but it is  unclear whether the wave is orbitally stable or not. 
  \subsubsection{Stability for the beam equation}
  We   state the relevant results for the beam equation. In it,  we need to make some assumptions regarding the existence of such waves and the spectrum of the corresponding operator $\ch$, defined in \eqref{b:12}. We should note that the existence of ground states, in a sense to be made precise below,  was proven by Levandosky, \cite{Lev}. Note however that \underline{our results apply to both the 
  periodic and the whole line cases}, with appropriate boundary conditions in each case.  
  \begin{theorem}
  \label{theo:beam}
  Let $p\geq 3$ be an odd integer and $I\subset (-\sqrt{2}, \sqrt{2})$ be an open interval.  Assume that there exists a one parameter 
  family $\{\vp_c\}_{c\in I}$  of solutions to \eqref{b:6} (where $L\geq 0$ could be finite  so that 
  \begin{itemize}
  \item $\vp_c\in H^1$ and $c\to \|\vp'_c\|_{L^2}$ is a differentiable function on $I$. 
  \item The operator $\ch_c=\p_x^4+c^2\p_x^2+1-p\vp_c^{p-1}$ satisfies \eqref{A} and \eqref{B}. 
  \end{itemize}
  Then, the wave $\vp_c$ is linearly stable if and only if 
   $\p_c\|\vp'_c\|< 0$ and 
  $$
 |c|\geq  \f{\|\vp'_c\|}{-2 \p_c \|\vp'_c\|}. 
  $$
  That is, if $\p_c\|\vp'_c\|\geq 0$ or $\p_c\|\vp'_c\|< 0$, but 
  $|c|<\f{\|\vp'_c\|}{-2 \p_c \|\vp'_c\|}$, the wave $\vp_c$ is unstable.  
  \end{theorem}
  {\bf Remarks: }
  \begin{enumerate}
  \item Levandosky, \cite{Lev}  has shown (in the case of the whole line) that   a family of ``ground state'' solutions may be constructed via a variational 
   procedure as follows. More precisely, one first solves the variational problem 
  \begin{equation}
  \label{var}
\left|\begin{array}{l} 
J(z)=\int_{-\infty}^\infty (z^2_{xx} -c^2z_x^2+z^2)dx\to \min \\
\textup{subject to} \ \ I(z)=\int_{-\infty}^\infty z^{p+1}(x) dx=1.          
\end{array} \right.  
  \end{equation}
Then, if one sets $\vp_c:=J(z_c)^{\f{1}{p-1}} z_c$, where $z_c$ is the solution to the minimization problem, then $\vp_c$ solves precisely \eqref{b:6}. Note here that the restriction $|c|<\sqrt{2}$ assures that such a problem has a solution and in fact $J(z)\geq 0$. 
Similar approach works in the periodic case as well, \cite{DS}. 
\item It should be noted however, that the ground state solutions produced via \eqref{var} need not be the only way of producing the family $\{\vp_c\}$ and in fact it is not, \cite{DS}. One reason for that  is the apparent lack of uniqueness in solving this type of minimization problem. Second, 
one may have ``excited'' states (or local minima of \eqref{var}), which are 
also interesting solutions to consider. We have not checked out  their stability, but it is likely that some of them are stable. We also point out that the variational methods for  studying the (orbital) stability of ground states are inapplicable for these excited states, whereas Theorem \ref{theo:beam} gives an exhaustive answer. 
  \item   We have shown in Section \ref{sec:1.2.1} that all reasonable 
  families of solutions must have the property that $\ch_c$ has a negative eigenvalue and a zero eigenvalue, with the eigenvector $\vp_c'$. The simplicity of such eigenvalues is however difficult to prove theoretically, due in particular to the lack of explicit formulas for $\vp_c$. Numerically however, one can easily check that $\ch_c$ has simple negative and simple zero eigenvalue, even for most of excited states, \cite{DS}. In fact, in all numerical runs in \cite{DS}, $\ch_c$ never failed to  satisfy \eqref{A} and \eqref{B}. 
  \item The restriction that $p$ is odd is technical and it is likely to 
  be unnecessary. More precisely, the issue is only the differentiability of the function $x\to |\vp(x)|^{p-1}\vp(x)$. 
  \end{enumerate}
The paper is organized as follows.   In Section \ref{sec:3.1}, we first present the relevant basic  results of the Shkalikov's theory for quadratic pencils, after which we introduce our main tool, the function $\cg$, list its main properties and derive  the proof of our main result, Theorem \ref{theo:5}. This is done modulo 
the somewhat technical Lemma \ref{le:c1} (which establishes the linear stability at $\om=\om^*(H)$), whose proof is given in Section \ref{sec:4.4}.  To achieve this, we   study and establish the Laurent expansion of the function $\cg$ at zero, which may be of independent use in future investigations. 

In the remaining sections, we apply the main result to obtain sharp results for linear stability/instability to the various examples announced above. In 
Section  \ref{sec:bous} we consider the Boussinesq model and prove Theorem \ref{theo:bous}. In Section \ref{sec:KGZ}, we set up the Klein-Gordon-Zakharov system and find the set of speeds that yields stable traveling waves as stated in Theorem \ref{theo:KGZ}.  Finally, in Section \ref{sec:beam}, we establish Theorem \ref{theo:beam} regarding the stability of the traveling waves for the 
 beam equation.

{\bf Acknowledgements:} We would like to take the opportunity to thank our colleagues and collaborators Sevdzhan Hakkaev, Panos Kevrekidis and Dmitry Pelinovsky for the numerous discussion on related topics.

\section{Proof of  Theorem \ref{theo:5}}
We first start with a preliminary material, which is the Shkalikov theory for quadratic pencils, \cite{Shkal}. 
\subsection{Some aspects of the Shkalikov's theory for quadratic pencils} 
\label{sec:3.1}
We follow the presentation of Shkalikov, \cite{Shkal}. We are especially interested in his index formula, which relates the number of unstable eigenvalues of $H$ with the  number of unstable eigenvalues of the pencil defined in \eqref{b:15}. 

Following Shkalikov, \cite{Shkal}, introduce a quadratic operator pencil in the form 
$$
A(\la)=\la^2 F+(D+i G) \la+T,
$$
where the coefficients $F,D,G,T$ are operators on a Hilbert space $H$, satisfying the following conditions 
\begin{enumerate}
\item[(i)] $F$ - bounded, invertible and self-adjoint
\item[(ii)] $(T,Dom(T))$ - self-adjoint, invertible
\item[(iii)] $D\geq 0,G$ are symmetric; $Dom(D), Dom(G)\subset Dom(T)$ and $D,G$ are $T-$bounded operators. 
\end{enumerate}
We say that an operator $M$ is $T-$bounded, if 
$$
|T|^{-1/2} M |T|^{-1/2}\in \cb(H)
$$
 We say that $\xi\in \rho(A)$, if $A(\xi)$, with domain $Dom(T)$ is invertible.

 Clearly, in our application, we will be interested in the case $F=Id$, $D=0$, 
 $G=-i\om \p_x$, $T=H$. However, {\bf note that $T=H$ is not invertible in our case of interest}. 
   Nevertheless, we introduce the associated quadratic pencil 
   $$
   \hat{A}(\la):=\la^2\hat{F}+\la(\hat{D}+i \hat{G})+J
   $$
 where $\hat{F}=|T|^{-1/2} F |T|^{-1/2}$, $\hat{G}=|T|^{-1/2} G |T|^{-1/2}$. We  will be also interested in the spectrum with respect to a smoother space. Namely, introduce the  space $H_{-1}$  with a norm 
 $$
 \|x\|_{H_{-1}}:=\| |T|^{1/2} x\|_{H}
 $$
It is shown in \cite{Shkal} that the spectrum of $A$ in $H_{-1}$ coincides with the spectrum of $\hat{A}$ considered on the space $H$. Note that in Definition \ref{defi:1}, we only consider smooth enough solutions $\psi$ anyway. Thus, we need to count the number of ``unstable'' eigenvalues of $\hat{A}$. 
 
The main result in the work of Shkalikov is Theorem 3.7. The following statement is a corollary of it. Here we have just presented a weaker version of the result, which will suffice for our purposes. 
\begin{theorem}(Theorem 3.7, \cite{Shkal})
\label{theo:shk}
Suppose the coefficients of the pencil $A$, $F,D,G,T$ satisfy conditions $(i)-(iii)$ listed above. Let the numbers of negative eigenvalues of $F$ and $T$, $\nu(F), \nu(T)$ respectively is finite. 

Then, the spectrum of $A(\la)$ in the open right-half plane $\cc_r=\{z: \Re z>0\}$, considered upon the space $H_{-1}$ consists of normal eigenvalues only. Moreover, the total algebraic multiplicity of all eigenvalues lying in $\cc_r$ satisfies 
$$
k(\hat{A})\leq \nu(T)+\nu(F). 
$$
\end{theorem}

 Recall that in our case however, the operator $T=H$ is not invertible. This case is also covered by Shkalikov, see Theorem 4.2,\cite{Shkal} under the structural assumption \eqref{A}. Our proof below just retraces his argument. 

Indeed, since $Ker(H)\neq \{0\}$, 
one needs to consider $H_\tau:=H+\tau Id$ for $0<\tau<<1$, so that $Ker(H_\tau)=\{0\}$.  

In order to apply Theorem \ref{theo:shk}, 
we need to check that $F=Id, G=-2 i \om \p_x$ are $H_\tau$ bounded. This amounts to showing that $|H_\tau|^{-1/2}\p_x |H_\tau|^{-1/2}\in \cb(L^2)$. This is a direct consequence of \eqref{E}. 

We can now apply Theorem \ref{theo:shk} to $\hat{A_\tau}$ to conclude 
$$
k(\hat{A_\tau})\leq \nu(H_\tau)+\nu(Id)=1,
$$
 for all small enough $\tau>0$.  
Since the eigenvalues depend continuously on $\tau$ (see \cite{Kato}, Chapter 7), we take a limit as $\tau\to 0+$ to get the desired inequality $k(\hat{A})\leq 1$. 

To recapitulate, we have shown that for fixed real $\om$ and under \eqref{A}, the equation 
\begin{equation}
\label{m:130}
\la^2\psi+2\la \om \psi'+H\psi=0, \ \ (\la, \psi)\in \rone_+\times L^2,
\end{equation}
has at most one solution with $\Re\la>0$. 

Finally, it is easy to see that in this case, the solution $\la$, if it exists, must be real. Indeed, suppose $\la\neq \bar{\la}$ is an eigenvalue for the pencil. That is, there is $\psi\in D(H):  \la^2\psi+2\la \om \psi'+H\psi=0$. Taking a complex conjugate and taking into account \eqref{G}, we see that 
$$
\bar{\la}^2\bar{\psi}+2\bar{\la} \om \bar{\psi}'+H \bar{\psi}=0
$$
Thus, $(\bar{\la}, \bar{\psi})$ is another solution to \eqref{m:130},with $\Re\bar{\la}>0$ and hence $k(A)\geq 2$, in contradiction with the inequality $k(A)\leq 1$. Thus, $\la$ must be real and we have established  
 \begin{corollary}
\label{cor:1} 
 For $\om\in \rone$, the equation  $\la^2 \psi+ 2\om  \la\psi'+H\psi=0$, $\la\in \cc, \psi\in D(H)$ 
 has at most one   solution   $(\la, \psi)$ with $\Re \la>0$. Moreover, such a pair will have $\la$ real, $\la>0$. 
  \end{corollary}

\subsection{Detecting instabilities: Proof of \eqref{m:2}} 
\label{sec:2.2}
Clearly, if  $\om=0$, then  \eqref{g1} is unstable, in the sense of Definition \ref{defi:1}, as in this case $\la=\de>0$ and $\psi=\phi$ will be a solution to \eqref{b:15}. So, assume henceforth that $\om \neq 0$. 

We start with an useful Proposition, which appears as Theorem in \cite{BO}, 
but it may in fact be a well-known result. 
\begin{proposition}
\label{BO} Assume that $A$ is a closed, densely defined (not necessarily self-adjoint) operator on a Hilbert space, which is bounded from below 
($\inf_{u\in D(A): \|u\|=1} \dpr{A u}{u}>-\infty$). Define its self-adjoint part $H=\Re A=\f{1}{2}[A+A^*]$. Then 
$$
\inf\{\Re\la: \la\in \si(A)\}\geq \inf \si(H). 
$$
In particular, if $H>0$, then $A$ is invertible. 
\end{proposition}

In  \eqref{g1}, we use the ansatz $u(t,x)=e^{\la t} \psi(x)$, whence one gets the following equation for $\psi$ 
\begin{equation}
\label{m:1}
\la^2 \psi+2\om \la \psi'+ H\psi=0
\end{equation}
Our first observation is that any nontrivial solution to \eqref{m:1} will have a non-trivial projection onto $\phi$. Indeed, assuming that $\psi\perp \phi$ and taking into account that $\{\phi\}^\perp$ is invariant under $H$, we can apply the projection operator $P_0$ on \eqref{m:1}. We get 
\begin{equation}
\label{psi}
(H+\la^2 +2\om \la P_0 \p_x P_0) \psi=0
\end{equation}
Considering now $(H+\la^2 +2\om \la P_0 \p_x P_0):\{\phi\}^\perp\to \{\phi\}^\perp$, we see that $(H+\la^2)|_{\{\phi\}^\perp}\geq \la^2$, whereas $2\om \la P_0 \p_x P_0$ is skew-adjoint operator on $\{\phi\}^\perp$. It follows that $\Re ((H+\la^2 +2\om \la P_0 \p_x P_0)=(H+\la^2)|_{\{\phi\}^\perp}\geq \la^2 ID$, whence $(H+\la^2 +2\om \la P_0 \p_x P_0)$ is invertible, by Proposition \ref{BO}. It follows from \eqref{psi} that $\psi=0$. Thus, $\dpr{\psi}{\phi}\neq 0$ for any non-trivial solution of \eqref{m:1}. 
  
Since solutions of \eqref{m:1} are up to multiplicative constant, we may take $\psi: \dpr{\psi}{\phi}=1$. That is  $\psi:=\phi+v$, where $v\perp\phi$.  We get the following equation for $v$ 
\begin{equation}
\label{m:5}
 (\la^2+2\om \la \p_x +H)v= (\de^2-\la^2)\phi-2\om \la \phi'
\end{equation}
Taking dot product with $\phi$ (and taking into account  $\dpr{\phi'}{\phi}=0$ and 
$\dpr{(H+\la^2)v}{\phi}=0$)  yields 
\begin{eqnarray*}
& & 2\om \la\dpr{v'}{\phi}=\de^2-\la^2 \\
& & \dpr{v}{\phi'}=\f{\la^2-\de^2}{2\om \la}
\end{eqnarray*}
Taking $P_0$ in \eqref{m:5} (and observing that $P_0 v=v$) on the other hand implies 
$$
[H+\la^2+2\om \la P_0\p_x P_0]v=-2\om \la \phi'.
$$
We now observe that the operator $\cl_\la=H+\la^2+2\om\la P_0\p_x P_0:\{\phi\}^\perp\to\{\phi\}^\perp $ is invertible by Proposition \ref{BO}  and the representation 
$$
\cl_\la=(H+\la^2)+\textup{skewsymmetric}:\{\phi\}^\perp\to\{\phi\}^\perp, \ \ H+\la^2|_{\{\phi\}^\perp}\geq \la^2 Id>0.
$$
Thus, $v=-2\om\la [H+\la^2+2\om\la P_0\p_x P_0]^{-1}[\phi']\in \{\phi\}^\perp$, since $\phi'\in \{\phi\}^\perp$. 
 
 From this analysis, we can say that the spectral problem \eqref{m:1} has a real solution $\la$   \underline{if and only if}
 \begin{equation}
 \label{m:10}
 \cg(\om;\la):=\dpr{[H+\la^2+2\om\la P_0\p_x P_0]^{-1}[\phi']}{\phi'}+\f{\la^2-\de^2}{4\om^2\la^2}
 \end{equation}
 has a positive root\footnote{Often times, we will omit the dependence of $\cg$ on  $\om$, since $\om$ will be mostly fixed} $\la_0>0$. In view of Corollary \ref{cor:1}, we have that this condition is necessary and sufficient for instability of \eqref{g1}. We have proved the following 
 \begin{proposition}
 \label{prop:shk}
 If $H$ satisfies the condition of theorem \ref{theo:5} and $\cg$ is as in \eqref{m:10}, then a  {\bf necessary and sufficient condition for 
 instability is that the function $\cg$ vanishes for some $\la_0>0$.}
 \end{proposition}
 Our next result claims   the continuity of the function $\cg(\om, \la)$ 
\begin{lemma}
\label{cont} Let $\{(\om_n, \la_n)\}\in \rone_+\times \rone_+$, so that 
$(\om_n, \la_n)\to (\om_0, \la_0)\in \rone_+\times \rone_+$. Then, 
$$
\lim_n \cg(\om_n; \la_n)= \cg(\om_0; \la_0)
$$
 \end{lemma}
 We provide the straightforward, but somewhat technical proof of 
 Lemma \ref{cont} in the Appendix. \\
 We now continue with a verification of the fact that $\cg(\om, \cdot)$ changes sign in $(0,\infty)$  (and hence vanishes somewhere there), provided either 
 $-\de^2\dpr{H^{-1}[\psi_0']}{\psi_0'}=|\dpr{\phi}{\psi_0'}|^2-\de^2 \|H^{-1/2} P_0[\psi_0']\|^2<0$ or $0\leq |\om|<\om^*(H)$. 
 \subsubsection{Analysis close to $\la=\infty$}
 We will show that under appropriate assumptions on $H$, we have 
 \begin{equation}
 \label{m:11}
 \lim_{\la\to \infty} \cg(\la)=\f{1}{4\om^2}>0.
   \end{equation}
 Indeed, $\lim_{\la\to \infty} \f{\la^2-\de^2}{4\om^2\la^2}= \f{1}{4\om^2}$, so it remains  
 $\lim_{\la\to \infty} \dpr{[H+\la^2+2\om\la P_0\p_x P_0]^{-1}[\phi']}{\phi'}=0$. 
 
 This requires a bit of elementary analysis of the operator $\cl$ defined above, in the limit of $\la\to \infty$. We have already established that $\cl_\la$ is invertible on $\{\phi\}^\perp$. We now need an estimate on the norm of its inverse. 
 \begin{proposition}
\label{prop:5}
For every $\la>0$, we have $H+\la^2+2\om \la P_0\p_x P_0:\{\phi\}^\perp\to \{\phi\}^\perp$. This operator has an inverse (in the said co-dimension one subspace) and its inverse satisfies   the estimate 
  \begin{equation}
  \label{m:75} 
  \|(H+\la^2+2\om\la P_0\p_x P_0)^{-1}\|_{\{\phi\}^\perp\to \{\phi\}^\perp} \leq \la^{-2}
  \end{equation}
\end{proposition}  
 \begin{proof}
 We have already checked the invertibility. Let  $g\in \{\phi\}^\perp$ is real-valued arbitrary function and 
 $f=(H+\la^2+2\om\la P_0\p_x P_0)^{-1}[g] \in \{\phi\}^\perp$, (note $f$ is real valued as well), so that 
  $$
  (H+\la^2+2\om\la P_0\p_x P_0) f=g,
  $$
 Taking dot product with $f$ yields  
  $$
  \la^2\|f\|^2\leq \dpr{(H+\la^2)f}{f}=\dpr{(H+\la^2+2\om\la P_0\p_x P_0) f}{f}=\dpr{f}{g} \leq \|f\| \|g\|,
  $$
  whence $\|(H+\la^2+2\om\la P_0\p_x P_0)^{-1}g\|=\|f\|\leq \la^{-2} \|g\|$, as stated. 
 \end{proof}

  Going back to the formula for $\cl_\la^{-1}$ (which now makes sense), we observe that   
 $$
 \limsup_{\la\to \infty}  |\dpr{[H+\la^2+2\om\la P_0\p_x P_0]^{-1}[\phi']}{\phi'}|\leq  \|\phi'\|^2 \limsup_{\la\to \infty} \la^{-2}  =0
 $$
 whence we have established \eqref{m:11}. 
 \subsubsection{Analysis close to $\la=0$}
 Regarding the behavior close to $\la=\ve\sim 0$, we observe first that 
 $$
 \f{\ve^2-\de^2}{4\om^2\ve^2}=-\f{\de^2}{4\om^2 \ve^2}+O(1),
 $$
 We will show that 
 \begin{equation}
 \label{m:15}
  \dpr{[H+\la^2+2\om\la P_0\p_x P_0]^{-1}[\phi']}{\phi'}=\f{1}{\ve^2} \f{\dpr{\phi'}{\psi_0}^2}{1+4\om^2\|H^{-1/2} P_0[\psi_0']\|^2}+O(\ve^{-1}).
   \end{equation}
 Before we prove \eqref{m:15}, let us pause for the moment to show that it implies all our statements for instabilities. First, if $\ 
 \de^{-2}|\dpr{\phi}{\psi_0'}|^2- \|H^{-1/2} P_0[\psi_0']\|^2=-\dpr{H^{-1}[\psi_0']}{\psi_0'}\leq 0$, 
   then the coefficient in front of $\ve^{-2}$ for $\cg(\ve)$ becomes 
   $$
   \f{\dpr{\phi'}{\psi_0}^2}{1+4\om^2\|H^{-1/2} P_0[\psi_0']\|^2}-\f{\de^2}{4\om^2}=
   \f{4\om^2(\dpr{\phi}{\psi'_0}^2-\de^2 \|H^{-1/2} P_0[\psi_0']\|^2)-\de^2}{4\om^2(1+4\om^2\|H^{-1/2} P_0[\psi_0']\|^2)}<0,
   $$
   whence $\lim_{\ve\to 0+}\cg(\ve)=-\infty$ and hence $\cg$ changes sign in $(0,\infty)$ and it has a root $\la_0>0$, hence instability. 
   
   If, on the other hand one has$ 
   |\dpr{\phi}{\psi_0'}|^2-\de^2 \|H^{-1/2} P_0[\psi_0']\|^2=-\de^2\dpr{H^{-1}[\psi_0']}{\psi_0'}>0$, the solution to the inequality 
   $$
   \f{\dpr{\phi'}{\psi_0}^2}{1+4\om^2\|H^{-1/2} P_0[\psi_0']\|^2}-\f{\de^2}{4\om^2}<0,
   $$
  in terms of $\om$ reads  as $0<|\om|<\om^*(H)$, that is \eqref{m:2}. But we observed that this implies change of sign for $\cg$. Thus there is $\la_0>0$, so that $\cg(\la_0)=0$, hence instability.

   Back to the proof of \eqref{m:15},  we need to compute the leading order term in 
   $$
   \dpr{(H+\ve^2+2\om\ve P_0\p_xP_0)^{-1}[\phi']}{\phi'}. 
   $$
Let $z\in \{\phi\}^\perp$ solves $(H+\ve^2+2\om\ve P_0\p_x P_0) z=\phi'$.  
Decompose $z=a\psi_0+  q$, $q\in \{\phi, \psi_0\}^\perp$.  
Observe that $q=P_1 q$. 
 The equation for $z$ becomes 
\begin{equation}
\label{m:25}
a\ve^2 \psi_0+2 a \om\ve P_0[\psi'_0]+(H+\ve^2+2\om\ve P_0\p_x P_0)q=\phi'
\end{equation}
The equation again has two different components - 
along $\psi_0$ and perpendicular to it. 

Along the direction of $\psi_0$, we take dot product with $\psi_0$ to obtain 
$$
a\ve^2+\dpr{(H+\ve^2+2\om \ve P_0\p_x P_0)q}{\psi_0}=\dpr{\phi'}{\psi_0}
$$
which implies (note $\dpr{q}{\psi_0}=0$, since $q=P_1 q\perp \psi_0$) 
\begin{equation}
 \label{m:20}
  a\ve^2 -2\om\ve  \dpr{ q}{\psi'_0}=\dpr{\phi'}{\psi_0}
   \end{equation}
   Along the orthogonal direction, we take $P_1$ in \eqref{m:25} - note that $P_1 P_0=P_1$ and $P_1 P_0 [\psi'_0]=P_0 P_1[\psi'_0]=P_0[\psi'_0]$). We get 
$$
 (H+\ve^2+2\om\ve P_1\p_x P_1)q= P_1[\phi']- 2 a \om \ve P_0[\psi'_0]
$$
Denoting the self-adjoint operator $H_1=H_1 P_1=P_1 H P_1\geq \si^2$ and 
consider $\cl_\ve= H_1+\ve^2 +2\om \ve P_1\p_x P_1$. Since $H_1+\ve^2\geq \si^2$, we may write 
$$
\cl_\ve=H_1[Id+ H_1^{-1}(\ve^2+2\om \ve P_1\p_x P_1)],
$$
whence by \eqref{E},  $\lim_{\ve\to 0+} \|H_1^{-1}(\ve^2+2\om\ve P_1\p_x P_1)\|_{L^2\to L^2}=0$, we can construct $\cl_\ve^{-1}$ in terms of Neumann series.  In fact 
$$
\cl_\ve^{-1}=H_1^{-1}+O(\ve).
$$
Since $H_1^{-1}=O(1)$, it follows that 
\begin{equation}
\label{m:30}
q=\cl_\ve^{-1}[- 2 a \om \ve P_0[\psi'_0]+P_1[\phi']]=-2 a \om \ve H_1^{-1}P_0[\psi'_0]+O(1)+a \ve O(\ve).
\end{equation}
Expressing $q$ from \eqref{m:30} back in \eqref{m:20} yields 
\begin{eqnarray*}
& & a\ve^2 -
2 \om \ve (-2 a \om \ve \dpr{H_1^{-1}P_0[\psi'_0]}{\psi'_0})+O(\ve)(1+a\ve^2)=\dpr{\phi'}{\psi_0}
\end{eqnarray*}   
Note 
$$
\dpr{H_1^{-1}P_0[\psi'_0]}{\psi'_0}=\dpr{H_1^{-1}P_0[\psi'_0]}{P_0[\psi'_0]}=\|H_1^{-1/2} 
[P_0[\psi'_0]]\|_{L^2}^2.
$$
As a consequence,  
$$
a\ve^2(1+4\om^2 \|H_1^{-1/2} [P_0[\psi'_0]]\|_{L^2}^2+O(\ve))=\dpr{\phi'}{\psi_0}+O(\ve),
$$
whence 
$$
a=\f{1}{\ve^2} \f{\dpr{\phi'}{\psi_0}}{1+4\om^2 \|H_1^{-1/2} [P_0[\psi'_0]]\|_{L^2}^2}+ O(\ve^{-1})
$$
Going back to \eqref{m:30}, we conclude that $q=O(\ve^{-1})$. Now 
\begin{eqnarray*}
 \dpr{(H+\ve^2+2\om\ve P_0\p_xP_0)^{-1}\phi'}{\phi'} &=& a\dpr{\psi_0}{\phi'}+\dpr{q}{\phi'}=\\
  &=& 
  \f{1}{\ve^2} \f{\dpr{\phi'}{\psi_0}^2}{1+4\om^2 \|H_1^{-1/2} [P_0[\psi'_0]]\|_{L^2}^2}+ O(\ve^{-1}). 
   \end{eqnarray*}   
   which establishes \eqref{m:15}. 
   
   We now progress to the proof of \eqref{m:3}. Introduce the notation 
 \begin{eqnarray*}
 \ca^{unstable} &:= & \{\om: \eqref{g1}\ \textup{is unstable}\}\\ 
 \ca^{stable} &:= & \{\om: \eqref{g1}\ \textup{is stable}\}\\
 \end{eqnarray*}
 We have already shown that $\ca^{unstable}\neq \emptyset$, in fact $(-\om^*(H), \om^*(H))\subset  \ca^{unstable}$. 
 
 The rest of the argument proceeds in two steps. First, we show that $(\om^*(H), \infty)\subset \ca^{stable}$ and then, we prove that $\om^*(H)\in \ca^{stable}$. 

 \subsection{Proof of $(-\infty, \om^*(H))\cup (\om^*(H), \infty)\subset \ca^{stable}$} 
\label{sec:3.3}

Let $\om_0: \om_0>  \om^*(H)>0$. 
As we saw, this condition implies that the function $\cg(\om_0; \la)$, defined in \eqref{m:10}, has 
$$
\lim_{\la\to \infty}\cg(\om_0; \la)=\f{1}{4\om_0^2}>0, \lim_{\la\to 0+} 
\cg(\om_0; \la)=+\infty.
$$
 By the continuity of the function $\la\to \cg(\om_0; \la)$, there are three  options. The first one is that $\cg(\om_0; \la)>0$ for all $\la>0$. By Proposition \ref{prop:shk}, this would mean that $\cg(\om_0, \la)$ never vanishes and hence, \eqref{m:1} does not have solutions, whence stability for $\om_0$ follows. We    claim that this is the only viable option, by disproving the other two. 
 
 The second option is that $\cg(\om_0; \la)$ has at least two different roots in $(0,\infty)$, say $0<\la_1<\la_2<\infty$. This implies that \eqref{m:1} will have two different solutions, $(\la_1, \psi_1), (\la_2, \psi_2)$, thus contradicting Corollary \ref{cor:1} of the Shkalikov index theory. We can therefore refute this option. 
 
 The third possibility  is that $\cg(\om_0;\la)$ has a double root, say $\la_0=\la_0(\om_0)>0$, but otherwise is a non-negative function. If we manage to 
 rule out this possibility, we will be done. This in fact not very hard. Assume for a contradiction that the double root actually occurs at $\la_0$. Let $\mu\neq 0$ and consider the associated perturbed problem for $H_\mu:=H-\mu^2\dpr{\cdot}{\phi}\phi$. Clearly, $P_0H_\mu P_0=H$ and thus,  according to \eqref{m:10} 
 \begin{equation}
 \label{67}
 \cg_{H_\mu}(\om, \la) = \dpr{[P_0 H_\mu P_0 +\la^2+2\om\la P_0\p_x P_0]^{-1}[\phi']}{\phi'}+\f{\la^2-\de^2-\mu^2}{4\om^2\la^2}=\cg(\om, \la)-\f{\mu^2}{4\om^2\la^2} 
 \end{equation}
The analysis that we performed for $H$ remains valid for $H_\mu$, at least for  
$0<\mu<<1$. That is, the root (if any) of $\la\to \cg_{H_\mu}(\om, \la)$ will give exactly the instability eigenvalues of  $H_\mu$. In particular, by the Shkalikov's theory, the number of zeros of $\cg_{H_\mu}(\om, \la)$ is either  one or zero.  

But if $\cg(\om_0, \la)$ has double zero at $\la_0$, it follows by 
\eqref{67} that for  all $\mu$ small enough, the function $\la\to \cd_{H_\mu}(\om_0, \la)$ will have at least two separate zeros $\la_1<\la_0< \la_2$, $\la_2-\la_1=O(\mu^2)$.  In any case, as we have alluded to above,  that is a contradiction with the Shkalikov's theory, which requires there is at most one unstable eigenvalue for $H_\mu$. With that, we have reached a contradiction, which implies that double roots are impossible for $\la\to \cg(\om, \la)$, when $\om>\om^*(H)$. Hence, 
$$
(\om^*(H), \infty)\subset \ca^{stable}.
$$

\subsection{Proof of $\om^*(H)\in \ca^{stable}$} 
\begin{lemma}
\label{le:c1}
The problem  \eqref{g1} is stable at $\om^*(H)$, i.e. $\om^*(H)\in \ca^{stable}$. 
\end{lemma}
 We now proceed to show Lemma \ref{le:c1}. Denote in this section $\om_0:=\om^*(H)$ for sake of brevity. We shall need several propositions. 
\begin{proposition}
\label{prop:15}
For any $\om$, the function $\ve\to \ve^2 \cg(\om,\ve)$ is real analytic (close to zero) and it has the Laurent expansion 
\begin{equation}
\label{h:1}
\cg(\om, \ve)=\ve^{-2} D_{-2}(\om)+\sum_{j=0}^\infty D_j(\om)\ve^j, 
\end{equation}
where $D_{-1}(\om)=0$, 
$$
D_{-2}(\om)= \f{\dpr{\phi'}{\psi_0}^2}{1+4\om^2\|H^{-1/2} P_0[\psi_0']\|^2}-\f{\de^2}{4\om^2},
$$
  and the functions $\{D_j\}^\infty_{j=1}$ are smooth functions of $\om$. Finally, the radius of analyticity $r(\om)$ satisfies 
  $$
  \inf_{\om\in [a,b]\subset(0, \infty)} r(\om)\geq r_{a,b}>0. 
  $$ 
  That is, on every compact interval $I=[a.b]$, one may choose a common radius of analyticity. 
\end{proposition}
\begin{proof}
Most of the statements in Proposition \ref{prop:15} were in fact considered before, see for example the formula for $D_{-2}(\om)$, derived in \eqref{m:15}. 

Let us first indicate the real analyticity of the function $\ve^2\cg(\ve)$. We have, from the defining equations and in our previous notations, 
$$
\cg(\om, \ve)=
a_\ve\dpr{\psi_0}{\phi'}+\dpr{q_\ve}{\phi'}+\f{1}{4\om^2}-\f{\de^2}{4\om^2\ve^2}. 
$$
where $a=a_\ve$ and $q=q_\ve$ are determined from the pair of equations 
\begin{eqnarray}
 \label{a}
& & a \ve^2-2\om \ve \dpr{q}{\psi_0'}=\dpr{\phi'}{\psi_0} \\
\label{q}
& & q_\ve=\cl_\ve^{-1}(-2a\ve\om P_0[\psi_0']+P_1[\phi']),
\end{eqnarray}
where $\cl_\ve=(H+\ve^2+2\om\ve P_1\p_x P_1)$. As we observed before, by the von Neumann expansion, 
$$
\cl_\ve^{-1}=\sum_{k=0}^\infty (-\ve)^k (H^{-1}(\ve+2\om P_1\p_x P_1))^k H_1^{-1}=
H_1^{-1}+O(\ve)
$$
and is real-analytic in $\ve$ in a region $|\ve|<\ve_0<<1$.  Substituting the expression for $q$ in the expression for $a$ yields the formula 
$$
a\ve^2+4\om^2\ve^2\dpr{\cl_\ve^{-1} P_0[\psi_0']}{[\psi_0']}=2\om\ve \dpr{\cl_\ve^{-1} P_1[\phi']}{[\psi_0']}+\dpr{\phi'}{\psi_0}
$$
and 
$$
a\ve^2= \f{2\om\ve \dpr{\cl_\ve^{-1} P_1[\phi']}{[\psi_0']}+\dpr{\phi'}{\psi_0}}{
1+4\om^2 \dpr{\cl_\ve^{-1} P_0[\psi_0']}{[\psi_0']}}.
$$
Observe that $\dpr{\cl_\ve^{-1} P_0[\psi_0']}{\psi_0'}=
\|H^{-1/2}P_0[\psi_0']\|^2+O(\ve)$, whence the denominator stays away from zero, whence  the real analyticity of $\ve\to \ve^2a(\ve)$. Furthermore, the coefficients $D_j(\om)$ (which can be written 
by means of  the Cauchy formula)  are smooth functions of $\om$.   
Substituting now this 
back in the formula \eqref{q} yields that $\ve\to \ve q(\ve)$ is real analytic too. Observe also that the radius of analyticity of both functions may be estimated below in terms of  $\om^{-1}$. Therefore, it is uniformly bounded away from zero, when $\om$ belongs to some compact interval, away from zero. 

We now compute $D_{-2}(\om), D_{-1}(\om)$. To that end, introduce 
$A, B$, so that 
$$
a(\ve)=\ve^{-2} A+\ve^{-1} B+O(1). 
$$
We compute $q$ up to order $O(1)$. We have $\cl_\ve^{-1}=H_1^{-1}+O(\ve)$ and therefore 
\begin{eqnarray*}
 q &=& (H_1^{-1}+O(\ve))(-2 \ve^{-1} A\om P_0[\psi_0']-2 B \om P_0[\psi_0']+P_1[\phi'])=\\
&=& -2 A \om \ve^{-1} H^{-1}P_0[\psi_0']-2B\om H^{-1} P_0[\psi_0']+H^{-1} P_1[\phi']+O(1)
\end{eqnarray*} 
Plugging this in \eqref{a}, we get two equations - at the scale of $\ve^0$ and $\ve^{1}$ 
\begin{eqnarray*}
& & A-2\om(-2 A \om\dpr{H_1^{-1}P_0[\psi_0']}{\psi_0'})=\dpr{\phi'}{\psi_0} \\
& & B-2\om(  \dpr{H^{-1} P_1[\phi']}{\psi_0'}-
2B\om\dpr{H^{-1} P_0[\psi_0']}{\psi_0'})=0. 
\end{eqnarray*}
From those, we get 
\begin{eqnarray*}
& & A = \f{\dpr{\phi'}{\psi_0}}{1+4\om^2 \|H^{-1/2}P_0[\psi_0']\|^2} \\
& & B= \f{2\om \dpr{H^{-1} P_1[\phi']}{\psi_0'}}{1+4\om^2 \|H^{-1/2}P_0[\psi_0']\|^2}
\end{eqnarray*}
We are now ready to compute $\cg(\om, \ve)$ up to order $O(1)$. We have 
\begin{eqnarray*}
\cg(\om, \ve) &=& 
a_\ve\dpr{\phi'}{\psi_0}+\dpr{q_\ve}{\phi'}+\f{1}{4\om^2}-\f{\de^2}{4\om^2\ve^2}= \\
&=& \ve^{-2}\left(A \dpr{\phi'}{\psi_0} -\f{\de^2}{4\om^2\ve^2}\right) 
+\ve^{-1}\left( B \dpr{\phi'}{\psi_0} -
2 A \om \dpr{H^{-1}P_0[\psi_0']}{\phi'}\right)+O(1)
\end{eqnarray*}
It remains to observe that 
\begin{eqnarray*}
D_{-2}(\om) &=& A \dpr{\phi'}{\psi_0} -\f{\de^2}{4\om^2}=\f{\dpr{\phi'}{\psi_0}^2}{1+4\om^2\|H^{-1/2} P_0[\psi_0']\|^2}-\f{\de^2}{4\om^2},\\
D_{-1}(\om) &=& B \dpr{\phi'}{\psi_0} -
2 A \om \dpr{H^{-1}P_0[\psi_0']}{\phi'}=0.
\end{eqnarray*}

\end{proof}
Several things to note. First, recall that the definition of $\om_0=\om^*(H)$ is so that $D_{-2}(\om_0)=0$.  In fact, $D_{-2}(\om)<0, \om<\om_0$ and 
$D_{-2}(\om)>0, \om>\om_0$. 

Now, consider the real-analytic function 
$$
\cg(\om_0,\ve)=\ve^{-2}D_{-2}(\om_0)+\sum_{j=0}^\infty D_j(\om_0)\ve^j= \sum_{j=0}^\infty D_j(\om_0)\ve^j. 
$$
By the analyticity, we have that there is $k$, so that $D_k(\om_0)\neq 0$. Indeed, otherwise, the real-analytic function  $\ve\to \cg(\om_0, \ve)$ will be identical to zero, at least for all $\ve$ small enough. This is however a contradiction with the meaning of the function $\cg(\om_0, \cdot)$. Indeed, not only it implies instability at $\om_0$, but it gives eigenvalues $\la=\ve\in (0, \ve_0)$, clearly in contradiction with our earlier conclusion (basically a consequence of Shkalikov theory) that there is at most one eigenvalue. 
Thus, there is $k: D_k(\om_0)\neq 0$. Let $k_0=\min\{k\geq 0:D_k(\om_0)\neq 0\}$.  \\
{\bf Claim:} $D_{k_0}(\om_0)>0$. \\
Indeed, take $\om_j\to \om_0+$. Since $\om_j\in \ca^{stable}$, take any $\ve: 0<\ve<r_0=\inf r(\om_j)$ (which is strictly positive by Proposition \ref{prop:15}), 
\begin{eqnarray*}
0\leq \limsup_j  \cg(\om_j, \ve) &=&  \limsup_j(\ve^{-2} D_{-2}(w_j)+(\sum_{m=0}^{k_0-1} D_m(\om_j)\ve^j)+D_{k_0}(\om_j)
\ve^{k_0}+O(\ve^{k_0+1}))= \\
&=& D_{k_0}(\om_0) \ve^{k_0}+O(\ve^{k_0+1})
\end{eqnarray*}
Clearly, the last inequality will be contradictory for small $\ve>0$, unless $D_{k_0}(\om_0)\geq 0$, which was the Claim (recall $D_{k_0}(\om_0)\neq 0$). 

Having established this claim, we easily finish the proof of our claim that $\om_0\in \ca^{stable}$. Indeed, we have $\lim_{\la\to \infty}\cg(\om_0, \la)=\f{1}{4\om_0^2}>0$ and 
$\lim_{\la\to 0+} \la^{-k_0} \cg(\om_0, \la)=D_{k_0}(\om_0)>0$. Again, assuming the contrary (i.e. $\om_0\in \ca^{unstable}$) implies the existence of $\la_0>0$, so that $\cg(\om_0, \la_0)=0$. By the behavior of $\la\to \cg(\om_0, \la)$ at $\la=0, \infty$ and the fact that the function $\cg(\om_0, \la)$ cannot have more than one zero in $(0, \infty)$, we again conclude that $\cg(\om_0, \cdot)$ has  double zero at $\la_0$. As we saw in Section \ref{sec:3.3}, this leads to a contradiction, if we consider the function 
$\la\to\cg_{H_\mu}(\om_0, \la)$. Namely, for all sufficiently small $\mu$, the function $\cg_{H_\mu}(\om_0, \la)$ will have at least two zeros, in contradiction with the Shkalikov's theory.

\subsection{Proof of Corollary \ref{cor:l1}}
In the course of our exposition for the Shkalikov's theory, we have proved that the spectrum of the pencil (and hence the spectrum of the operator $\ct$) is invariant under complex conjugation. More specifically, if $(\la, \psi)$ solves \eqref{g1}, then so does $(\bar{\la}, \bar{\psi})$. 

Now, under the assumptions in Corollary \ref{cor:l1}, we can also show that the spectrum is invariant under the transformation $\la\to -\la$. Indeed, if $(\la, \psi)$ solves \eqref{g1}, then so does $(-\la, \psi(-\cdot))$. 

Thus, in the unstable situation, we have shown that there is unique positive eigenvalue $\la$ of $\ct$ and hence there is an unique negative one as well. Moreover, this argument excludes the possibility for more negative spectrum, since this would imply more positive (unstable) eigenvalues for the pencil, which is forbidden by the Shkalikov's theory. 

In the stable situation, one shows by contradiction argument, similar to the one in the previous paragraph, that there is no spectrum off the imaginary axes and thus $\si(\ct)\subset i \rone$.

 \section{Linear stability analysis for the Boussinesq model} 
 \label{sec:bous}
  The proof of Theorem \ref{theo:bous} follows directly from Theorem \ref{theo:5}, applied to the linearized problem $z_{tt}+2c z_{t x}+H z=0$, where $H$ is as defined in \eqref{z:2}. In order to follow this program, we need to verify that the conditions in Theorem   \ref{theo:5} for the operator $H$ are met, then we need to compute the critical index $\om^*(H)$ and finally, we need to establish that the linear instability for $z_{tt}+2c z_{t x}+H z=0$ is equivalent to the linear instability for the actual problem \eqref{zz:2}. We will do this in a sequence of several propositions. 
  \begin{proposition}(Spectral properties of $H$)
  \label{pb:1} The operator $H$, defined in \eqref{z:2} satisfies \eqref{A}, \eqref{E}, \eqref{B}, \eqref{G}.    
  \end{proposition}
  Next, 
  \begin{proposition}
  \label{pb:2} The critical index $\om^*(H)=\f{\sqrt{(p-1)(1-c^2)}}{\sqrt{5-p}}.$ 
  
  \end{proposition}
  Regarding the existence of the unstable solutions, we have already explained that an instability in $z_{tt}+2c z_{t x}+H z=0$ implies instability in \eqref{zz:2}. We need to show then the reverse. 
\begin{proposition}
  \label{pb:3}  
  Suppose that $V\in H^4(\rone)$ solves the equation 
  $$
  \la^2 V+2 c \la V'+ T[V]=0,
  $$
  for some $\la>0$. Then there exists a function $Z\in H^4(\rone)$, so that $Z'=V$ and so that $Z$ solves 
  $$
  \la^2 Z+2c\la Z'+ H[Z]=0
  $$
  \end{proposition}
Let us show now that Theorem \ref{theo:bous} is a simple consequence of Propositions \ref{pb:1}, \ref{pb:2} and \ref{pb:3}. Indeed, since we can apply Theorem \ref{theo:5}, we only need to find the intervals, in which the speeds yield stable traveling waves. To that end, we set up the inequality 
$$
1>|c|\geq \om^*(H)=\f{\sqrt{(p-1)(1-c^2)}}{\sqrt{5-p}}. 
$$
The solution to this inequality is  
$$
1>|c|\geq \f{\sqrt{p-1}}{2}. 
$$
In the complementary set, $0\leq |c|<\f{\sqrt{p-1}}{2}$, 
 we have instability.  This finishes the proof of Theorem \ref{theo:bous}, modulo the claims of the Propositions.

\subsection{Proof of Proposition \ref{pb:1}} 

We start off by noting that \eqref{G} is obvious by inspection. Similarly, \eqref{E} follows from the fact that $H$ has one negative eigenvalue\footnote{to be proved below in this proposition} and thus, for $\tau>>1$, we have that $(H+\tau)>Id$, hence is invertible and moreover 
$$
(H+\tau)^{-1}:L^2\to H^4,
$$
and more generally for all $s\in [0,1]$, $(H+\tau)^{-s}:L^2\to H^{4s}$. Therefore, we conclude $(H+\tau)^{-1/2}\p_x (H+\tau)^{-1/2}, \p_x (H+\tau)^{-1} :L^2\to L^2$ are a bounded operators. 

\subsubsection{Proof of \eqref{A}} 
\label{sec:A}  To that end, note first that $H$ can be 
written in the following form 
$$
H=-\p_x\cl \p_x, \ \ \cl:=-\p_x^2+(1-c^2)-p\vp_c^{p-1} 
$$
The operator $\cl$ is a second order differential operator, which appears in the study of the KdV equation in exactly the same way - as a linearization around the traveling wave solution $\vp_c$. Its spectral properties are well-documented and we just remind the ones that are of concern to us.  Namely it has a simple eigenvalue at zero (with eigenvector $\vp_c'$), single and simple negative eigenvalue and it has a  spectral gap, that is the rest of the spectrum lies in a set $[\ka, \infty), \ka>0$, see for example \cite{CGNT}, Section 3.  
 
We will now be able to infer the same structure for the spectrum of $H$, mainly due to the representation $H=-\p_x\cl \p_x$. First, it has an eigenvalue at zero, with eigenvector $\vp_c$ ($H\vp_c=-\p_x \cl[\vp_c']=0$). This eigenvalue 
is simple (otherwise 
a contradiction with the simplicity of the zero eigenvalue for $\cl$). 

Next, we check that $H$ has a negative eigenvalue. Taking two spatial derivatives in \eqref{s:1} we get 
\begin{equation}
\label{s:5}
\vp_c''''-(1-c^2)\vp_c''+p \vp_c^{p-1} \vp''+p(p-1)\vp_c^{p-2} (\vp_c')^2=0
\end{equation}
taking dot product with $\vp_c''$ and taking into account that 
$-\cl \vp_c''=\vp_c''''-(1-c^2)\vp_c''+p \vp_c^{p-1} \vp''$, we get 
$$
-\dpr{\cl \vp_c''}{\vp_c''}+p(p-1)\int \vp_c^{p-2} (\vp_c')^2\vp_c'' dx=0
$$
But $\int \vp_c^{p-2} (\vp_c')^2\vp_c'' dx=-\f{p-2}{3} \int  \vp_c^{p-3} (\vp_c')^2dx$. Thus 
$$
\dpr{H \vp_c'}{\vp'_c}=\dpr{\cl \vp_c''}{\vp_c''}=-\f{p(p-1)(p-2)}{3} \int  \vp_c^{p-3} (\vp_c')^2dx.
$$
Note that if $p>2$, the expression is negative, thus yielding directly (by the Ritz-Reileigh criteria) that the minimal eigenvalue 
$$
\la_0(H)=\inf_{\psi:\|\psi\|=1} \dpr{H \psi}{\psi}\leq \|\vp_c'\|^{-2} \dpr{H \vp_c'}{\vp'_c}<0. 
$$
Even in the case $p=2$, this argument implies the existence of a negative eigenvalue. Indeed, assuming that $H\geq 0$ (and since we already know that $\vp_c$ is a simple eigenvalue), it follows that $H|_{\{\vp_c\}^\perp}>0$. Thus, since $\vp_c'\in \{\vp_c\}^\perp$, it should be that $\dpr{H \vp_c'}{\vp'_c}>0$, a contradiction with $\dpr{H \vp_c'}{\vp'_c}=0$, in the case $p=2$. This shows that $\la_0(H)<0$. 

Finally, we need to establish that $\la_0(H)$ is simple. This would follow, if we manage to show that $\la_1(H)\geq 0$ (and indeed, since $0$ is an eigenvalue, it would follow that $\la_1(H)=0$). Denote the (smooth) eigenvector of $\cl$ by $\zeta$. One may apply the Courant maxmin principle for the first eigenvalue, which states 
$$
\la_1(H)=\sup_{z\neq 0} \inf_{u\perp z} \f{\dpr{H u}{u}}{\|u\|^2}
$$
Taking $z=\phi'$ above yields 
$$
\la_1(H)=\sup_{z\neq 0} \inf_{u\perp z} \f{\dpr{H u}{u}}{\|u\|^2}\geq \inf_{u\perp \phi'} \f{\dpr{H u}{u}}{\|u\|^2}= \inf_{u\perp \phi'} \f{\dpr{\cl u'}{u'}}{\|u\|^2}=\inf_{u'\perp \phi} \f{\dpr{\cl u'}{u'}}{\|u\|^2}
$$
where in the last identity, we have used that $u\perp \phi'$ exactly when $u'\perp \phi$.  
It remains to observe now that since $\cl$ has a simple negative eigenvalue, with eigenvector $\phi$, we have $\cl|_{\{\phi\}^\perp}\geq 0$ and hence $\dpr{\cl u'}{u'}\geq 0$ (since $u'\in \{\phi\}^\perp$), whence 
$$
\la_1(H)\geq \inf_{u'\perp \phi} \f{\dpr{\cl u'}{u'}}{\|u\|^2}\geq 0. 
$$ 
Thus, property \eqref{A} is fully established. 
\subsubsection{Proof of \eqref{B}} 
Note first that in our notations, $\psi_0=\vp_c/\|\vp_c\|$ and denote the negative eigenvector of $H$ by $\phi, \|\phi\|=1$. We need to show $\dpr{\phi'}{\vp_c}\neq 0$. 
We separate the proof in the cases $p>2$ and $p=2$. We have by our computations in Section \ref{sec:A} above that 
$$
 \dpr{H \vp_c'}{\vp_c'}<0
$$
It follows that  
$$
0>\dpr{H \vp_c'}{\vp_c'}=\la_0 \dpr{\vp'_c}{\phi}^2+
\dpr{H( \vp_c'-\dpr{\vp_c'}{\phi}\phi)}{(\vp_c'-\dpr{\vp_c'}{\phi}\phi)} 
$$
Note that $(\vp_c'-\dpr{\vp_c'}{\phi}\phi) \in \{\phi\}^\perp$. Since $\la_0(H)$ is the only negative eigenvalue for $H$ and it is simple, it follows that $H|_{\{\phi\}^\perp}\geq 0$ and in particular 
$$
\dpr{H( \vp_c'-\dpr{\vp_c'}{\phi}\phi)}{(\vp_c'-\dpr{\vp_c'}{\phi}\phi)}\geq 0. 
$$
The last two inequalities imply that 
$$
0>\la_0 \dpr{\vp'_c}{\phi}^2,
$$
whence $\dpr{\vp'_c}{\phi}\neq 0$. But then $\dpr{\phi'}{\vp_c}=-\dpr{\vp'_c}{\phi}\neq 0$. 

In the case $p=2$, as we have shown $\dpr{H \vp_c'}{\vp_c'}=0$, so we need to be more precise in the arguments above. In particular, we need to observe that since $\vp_c'\perp \vp_c$ and $\phi\perp\vp_c$ (as eigenvectors corresponding to different eigenvalues), it follows that actually 
\begin{equation}
\label{s:10}
\vp_c'-\dpr{\vp_c'}{\phi}\phi \in span[\phi, \vp_c]^\perp.
\end{equation}
Moreover, by \eqref{s:5} (in the case $p=2$),  
we have that $\cl[\vp_c'']=2(\vp_c')^2$, whence 
$$
\ch[\vp_c']=-\p_x \cl[\vp_c'']=-4\vp_c'\vp_c''\neq \la_0 \vp_c'. 
$$ 
The last computation shows that $\vp_c'$ and $\phi$ are linearly independent and hence, we can upgrade \eqref{s:10} to 
\begin{equation}
\label{s:12}
0\neq \vp_c'-\dpr{\vp_c'}{\phi}\phi \in span[\phi, \vp_c]^\perp.
\end{equation}
By \eqref{A}, in particular the spectral gap that we have established for $H$,  and \eqref{s:12},   it follows that 
$$
\dpr{H(\vp_c'-\dpr{\vp_c'}{\phi}\phi)}{(\vp_c'-\dpr{\vp_c'}{\phi}\phi)}\geq 
\ka  \|\vp_c'-\dpr{\vp_c'}{\phi}\phi\|^2>0,
$$
where $\ka$ is the size of the spectral gap.  We can conclude 
\begin{eqnarray*}
\la_0 \dpr{\vp'_c}{\phi}^2 &=& \dpr{H \vp_c'}{\vp_c'}-\dpr{H(\vp_c'-\dpr{\vp_c'}{\phi}\phi)}{(\vp_c'-\dpr{\vp_c'}{\phi}\phi)}=\\
&=& -\dpr{H(\vp_c'-\dpr{\vp_c'}{\phi}\phi)}{(\vp_c'-\dpr{\vp_c'}{\phi}\phi)}<0,
\end{eqnarray*}
whence $\dpr{\phi'}{\vp_c}=-\dpr{\vp'_c}{\phi}\neq 0$. 

\subsection{Proof of Proposition \ref{pb:2}} 
Recall $H=-\p_x\cl\p_x$. In our notations, we need to compute 
$\dpr{H^{-1} \psi_0'}{\psi_0'}=\f{1}{\|\vp_c\|^2} \dpr{H^{-1} \vp_c'}{\vp_c'}$.

Again, starting from \eqref{s:1}, take a derivative in the parameter $c$. We obtain, 
$$
(\p_c\vp)''-(1-c^2)\p_c \vp+p\vp^{p-1}\p_c \vp+2c\vp=0.
$$
In terms of the operator $\cl$, we have $\cl[\p_c \vp]= 2c\vp$ or 
$\f{1}{2c}\p_c \vp=  \cl^{-1}[\vp_c]$. Here recall that $\cl$ is invertible on $Ker[\cl]^\perp=\{\vp_c'\}^\perp$ and $\vp_c\in \{\vp_c'\}^\perp$. 

Introduce $z$, so that $\vp_c'=H[z]$. We know that such a $z$ exists, since the operator $H$ is invertible on the subspace  $Ker[H]^\perp=\{\vp_c\}^\perp$ and $\vp_c'\in \{\vp_c\}^\perp$. We have 
$$
\vp_c'=H[z]=-\p_x\cl[z']
$$ 
Thus $\vp_c=-\cl[z']+const$ and the constant turns out to be zero by testing this identity at $x\to \infty$. Thus, $z$ is such that $\cl[z']=-\vp_c$ and 
hence 
$$
z'=-\cl^{-1}[\vp_c]=-\f{1}{2c}\p_c \vp. 
$$
Now 
$$
-\dpr{H^{-1} \vp_c'}{\vp_c'}=-\dpr{z}{\vp_c'}=\dpr{z'}{\vp_c}=
-\f{1}{2c}\dpr{\p_c \vp_c}{\vp_c}=-\f{1}{4c}\p_c \|\vp_c\|^2. 
$$
All in all 
$$
-  \dpr{H^{-1} \psi_0'}{\psi_0'}=
-\f{1}{\|\vp_c\|^2} \dpr{H^{-1} \vp_c'}{\vp_c'}=-\f{1}{4c} \f{\p_c[\|\vp_c\|^2]}{\|\vp_c\|^2]}
$$
It remains to compute $\|\vp_c\|^2$ and perform the elementary calculus operations. We have 
\begin{eqnarray*}
\|\vp_c\|^2 &=& \left(\f{p+1}{2}(1-c^2)\right)^{\f{2}{p-1}}
\int_{-\infty}^\infty sech^{\f{4}{p-1}}\left(\f{\sqrt{1-c^2}(p-1)}{2}\xi\right) d\xi= \\
&=& m_p (1-c^2)^{\f{2}{p-1}-\f{1}{2}} 
\end{eqnarray*}
Thus, 
$$
-\f{1}{4c} \f{\p_c[\|\vp_c\|^2]}{\|\vp_c\|^2]}= \f{5-p}{4(p-1)(1-c^2)}. 
$$
Clearly, for $p\geq 5$, we have $\dpr{H^{-1} \psi_0'}{\psi_0'}=\f{p-5}{4(p-1)(1-c^2)}\geq 0$ and hence instability. 
If $p<5$,  we get 
$$
\om^*(H)= \f{1}{2\sqrt{-  \dpr{H^{-1} \psi_0'}{\psi_0'}}}= 
\f{\sqrt{(p-1)(1-c^2)}}{\sqrt{5-p}}.
$$

\subsection{Proof of Proposition \ref{pb:3}} 
Proposition \ref{pb:3} is a simple exercise. Writing down the particular form of the operator $T$, we see that $V$ satisfy 
$$
(\p_x^4-(1-c^2)\p_x^2+2c\la\p_x +\la^2)V=-p \p_x^2(\vp_c^{p-1} V). 
$$
Observe that the symbol of the operator $M=\p_x^4-(1-c^2)\p_x^2+2c\la\p_x +\la^2$ 
$$
\xi^4+(1-c^2)\xi^2-2c i\la\xi+\la^2
$$
has a positive real part, bounded away from zero and is thus invertible on $L^2$.  In particular, we may write 
$$
V=-p M^{-1}\p_x^2 [\vp_c^{p-1} V].
$$
Now, it suffices to set 
$$
Z:=-p M^{-1}\p_x [\vp_c^{p-1} V],
$$
so that $Z_x=V$ (and hence it satisfies $\la^2 Z+2c\la Z+H[Z]=0$) and check that such a function is well-defined and belongs\footnote{further regularity of course holds and can be inferred in the standard bootstrap fashion} to $L^2(\rone)$. Since $\vp_c^{p-1}V\in L^2$, it suffices to show that $\p_x M^{-1}$ defines a bounded operator on $L^2$ or 
$$
\sup_{\xi} \left|\f{\xi}{\xi^4+(1-c^2)\xi^2-2c i\la\xi+\la^2}\right|<\infty
$$
The last inequality follows by inspection.

\section{Linear stability analysis for the Klein-Gordon-Zakharov system} 
\label{sec:KGZ}
In this section, we will provide the complete details of the proof of Theorem \ref{theo:KGZ}. We need to first setup the problem at hand in the form \eqref{b:15}. 
This is not straight forward, due to the requirement that $\ch$ be self-adjoint. In fact, the choice of $\f{1}{2}$ that we have made in the formulation of \eqref{KGZ} helps us accomplish exactly that.  

We take the ansatz $u(t,x)=\vp(x-ct)+v(t,x-ct), n(t,x)= \psi(x-ct)+w(t, x-ct)$ for real-valued functions $v,w$. We plug it in \eqref{KGZ} and ignore all terms $O(v^2+w^2)$. We get the following linear equation for $(v,w)$, 
\begin{equation}
\label{KGZ1}
\left\{
\begin{array}{l}
v_{tt}-2cv_{tx} -(1-c^2)v_{xx}+v+\psi v+ \vp w =0 \\
w_{tt}-2cw_{tx}-(1-c^2) w_{xx}-(\vp v)_{xx}=0. 
\end{array}
\right.
\end{equation}
At this point, we introduce the convenient quantity $\mu=\sqrt{1-c^2}$. 
Setting $w=z_x$ and taking off one derivative from the second equation, 
 we obtain the following system for $\vec{\Phi}=\left(\begin{array}{c} 
v \\
z
\end{array} \right)$
\begin{equation}
\label{KGZ2}
\vec{\Phi}_{tt}-2c\vec{\Phi}_{t x}+\ch \vec{\Phi}=0, \ \ \ch:=\left(\begin{array}{c c} 
H_1 & A \\
A^* & H_2
\end{array} \right), 
\end{equation}
where 
\begin{eqnarray*}
& & H_1 = -(1-c^2)\p_{xx}+1+\psi=-\mu^2\p_{xx}+1-\f{\vp^2}{2\mu^2}\\
& & H_2=-(1-c^2)\p_{xx}=-\mu^2\p_{xx}\\
& & Az=\vp z_x, A^* v=-(\vp v)_x
\end{eqnarray*}
Clearly, the operator $\ch$ is self-adjoint. It now remains to show that $\ch$ satisfies the requirements of Theorem \ref{theo:5}, after which, we will compute the quantity of interest $\dpr{\ch^{-1}\psi_0'}{\psi_0'}$. 
We formulate the needed results in a series of Propositions. 
\begin{proposition}
\label{prop:K1}
The self-adjoint operator $\ch$ defined in \eqref{KGZ2} has a simple eigenvalue at zero, with an eigenvector $(\vp', -\f{\vp^2}{2\mu^2})$. 
\end{proposition}

\begin{proposition}
\label{prop:K2}
The self-adjoint operator $\ch$ defined in \eqref{KGZ2} has one  simple negative eigenvalue. 
\end{proposition}

\subsection{Proof of Proposition \ref{prop:K1}}
We need to solve $\ch \left(\begin{array}{c} f \\ g \end{array}   \right)=0$. This is equivalent to 
$$
\left|
\begin{array}{l}
-\mu^2 f''+f -\f{\vp^2}{2 \mu^2} f+ \vp g'=0 \\
-(\vp f)'-\mu^2 g''=0.
\end{array}
\right.
$$
Integrating the second equation yields  $g'=-\f{\vp f}{\mu^2}$, which we plug in the first equation. We get 
\begin{equation}
\label{k2}
-\mu^2 f''+f-\f{3\vp^2}{2\mu^2} f=0.
\end{equation}
Recall now the equation \eqref{eq:k1}, which defines $\vp$. In fact, taking a derivative $\p_x$ in it yields 
$$
-\mu^2 \vp'''+\vp'-\f{3\vp^2}{2\mu^2} \vp'=0. 
$$
By comparing the last two formulas, we see that $f=\vp'$ is a solution of \eqref{k2}. Moreover, by the standard theory for one-dimensional Hill's operators\footnote{Note that this particular operator is in fact the ubiquitous $L_-$, which appears in the linearization of standing waves for the cubic Schr\"odinger equation, which is known to have one-dimensional kernel} (see for example \cite{CGNT}, Section 3), $\vp'$ is the unique solution to \eqref{k2}, up to a multiplicative constant. 

It remains to observe that since $g'=-\f{\vp f}{\mu^2}=-\f{\vp \vp'}{\mu^2}$, we have 
$$
g=-\f{\vp^2}{2\mu^2}
$$
and thus the eigenvector, corresponding to the zero eigenvalue is $(\vp', -\f{\vp^2}{2\mu^2})$. 

\subsection{Proof of Proposition \ref{prop:K2}}
The proof of Proposition \ref{prop:K2} is less standard than the proof of Proposition \ref{prop:K1}. We set up the eigenvalue problem in the form 
$\ch \left(\begin{array}{c} f \\ g \end{array}   \right)=-a^2 \left(\begin{array}{c} f \\ g \end{array}   \right)$ for some $a\in (0,\infty)$. That is, we need to show that there exists an unique $a_0>0$, so that the eigenvalue problem has an unique (up to a multiplicative constant) solution $\left(\begin{array}{c} f \\ g \end{array}   \right)$.  Write the eigenvalue problem in the form  
$$
\left|
\begin{array}{l}
-\mu^2 f''+f -\f{\vp^2}{2 \mu^2} f+ \vp g'=-a^2 f \\
-\mu^2 g''-(\vp f)' =-a^2 g.
\end{array}
\right.
$$
From the second equation, we express $g=\p_x(a^2-\mu^2\p_{xx})^{-1}[\vp f]$. 
This is possible, since \\ $a^2-\mu^2\p_{xx}\geq a^2 Id$ and hence is invertible. Thus, 
$$
g'=\p_{xx}(a^2-\mu^2\p_{xx})^{-1}[\vp f]=-\f{\vp f}{\mu^2}+\f{a^2}{\mu^2} (a^2-\mu^2\p_{xx})^{-1}[\vp f].
$$
Plugging this in the first equation, we obtain 
\begin{equation}
\label{k3}
-\mu^2 \p_{xx} f+(1+a^2) f- \f{3\vp^2}{2\mu^2} f+\f{a^2}{\mu^2}[\vp  (a^2-\mu^2\p_{xx})^{-1}(\vp f)]=0. 
\end{equation}
To recapitulate, we have obtained the equation \eqref{k3} to be equivalent to the eigenvalue problem $\ch \left(\begin{array}{c} f \\ g \end{array}   \right)=-a^2 \left(\begin{array}{c} f \\ g \end{array}   \right)$. That is, to prove Proposition  \ref{prop:K2}, we need to show that there exists an unique  $a_0>0$, so that the operator 
$$
M_a=-\mu^2 \p_{xx} +(1+a^2) - \f{3\vp^2}{2\mu^2} +\f{a^2}{\mu^2}[\vp  (a^2-\mu^2\p_{xx})^{-1}(\vp \cdot)],
$$
has an eigenvalue at zero for $a=a_0$ and in addition, this eigenvalue is simple.  

Several things to note for the one-parameter family of operators $M_a$. It is clear that $M_a, a\geq 0$ are  self-adjoint and in addition, \\
\\
{\bf Claim:} If $ a\geq b\geq 0$, then $M_a\geq M_b+(a^2-b^2)Id\geq M_b$. \\
\\ 
Let us finish the proof of Proposition \ref{prop:K2}, based on this Claim, whose proof we postpone for the end of this section. Denote 
$$
\la(a)=\inf\{\la: \la\in \si(M_a)\}=\inf_{\|f\|=1} \dpr{M_a f}{f}.
$$
Clearly, $a\to \la(a)$ is a continuous function and in view of Claim, $a\to \la(a)$ is an increasing function of its argument.  

Next, we consider the easy case $a>>1$. Observe that if $a^2>\sup_x \f{3\vp^2(x)}{2\mu^2}-1$, $M_a$ becomes a positive operator (the operator $[\vp  (a^2-\mu^2\p_{xx})^{-1}(\vp \cdot)]$ is positive by its Fourier transform representation, see the proof of the Claim below) and thus $\la(a)>0$. 
 
 The case $a=0$ presents 
another interesting observation. Namely, since 
$$
M_0=-\mu^2 \p_{xx} +1 - \f{3\vp^2}{2\mu^2},
$$
we already know that $M_0[\vp']=0$ and hence by Sturm-Liouville theory, there is a negative eigenvalue for $M_0$, that is $\la(0)<0$. Alternatively, one may directly compute 
$$
\dpr{M_0 \vp}{\vp}=-\f{1}{\mu^2}\int_{-\infty}^\infty \vp^3(y) dy<0,
$$
whence again $\la(0)<0$. Thus, the continuous and increasing 
function $a\to \la(a)$ is negative at $a=0$ and positive for large $a$, whence it has exactly one zero, say $a_0$. Thus, the eigenvalue that we are looking for is $\la=-a_0^2<0$. We still need to check that this eigenvalue is simple, which by the equivalences that we have established means that we have to show that $0$ is an isolated (simple) eigenvalue of $M_{a_0}$. 

This is however 
an easy consequence of the fact that $M_{a_0}\geq M_0+a_0^2 Id$ (from the Claim). Indeed, denote by $\phi_0$ the eigenvector for $M_0$, which corresponds to its unique and simple negative eigenvalue (\cite{CGNT}, Section 3). Then, 
 $M_0|_{\{\phi_0\}^\perp}\geq 0$, whence  by the Courant maxmin principle
$$
\la_1(M_{a_0})=\sup_{z\neq 0} \inf_{u\perp z: \|u\|=1} \dpr{M_{a_0} u}{u} \geq 
a_0^2+\inf_{u\perp \phi_0: \|u\|=1} \dpr{M_{0} u}{u}\geq a_0^2,
$$
Thus, we have  that $\la(a_0)=\la_0(M_{a_0})=0$, 
while $\la_1(M_{a_0})\geq a_0^2$, which finishes the proof of Proposition \ref{prop:K2}. 

\subsubsection{Proof of Claim}
By the particular form of the operators $M_a$, it suffices to establish that for all test functions $f$, 
\begin{equation}
\label{k5}
a^2 \dpr{[\vp  (a^2-\mu^2\p_{xx})^{-1}(\vp f)]}{f} \geq b^2 \dpr{[\vp  (b^2-\mu^2\p_{xx})^{-1}(\vp f)]}{f},
\end{equation} 
whenever $a\geq b$. This is easily seen by the Plancherel's theorem. 
Denoting $h=\vp f$, we have 
\begin{eqnarray*}
& & a^2 \dpr{[\vp  (a^2-\mu^2\p_{xx})^{-1}(\vp f)]}{f}=
\int_{-\infty}^\infty \f{a^2}{a^2+4\pi^2\mu^2 \xi^2} |\hat{h}(\xi)|^2 d\xi \geq \\
& & \geq \int_{-\infty}^\infty \f{b^2}{b^2+4\pi^2\mu^2 \xi^2} |\hat{h}(\xi)|^2 d\xi=
 b^2 \dpr{[\vp  (b^2-\mu^2\p_{xx})^{-1}(\vp f)]}{f},
\end{eqnarray*}
where we have used the elementary inequality $\f{a^2}{a^2+4\pi^2\mu^2 \xi^2}\geq \f{b^2}{b^2+4\pi^2\mu^2 \xi^2}$, if $a\geq b$.

\subsection{Proof of Theorem \ref{theo:KGZ}}
Now that we have checked most of the required conditions on $\ch$, let us verify the remaining ones and the quantity $\dpr{\ch^{-1}\psi_0'}{\psi_0'}$. Here  
$\psi_0=m\left(\begin{array}{c} \vp' \\ -\f{\vp^2}{2\mu^2} \end{array}\right)$, where $m$ is so that $\|\psi_0\|=1$. Thus, we need to compute $\ch^{-1}\left(\begin{array}{c} \vp' \\ -\f{\vp^2}{2\mu^2} \end{array}\right)$. Recall that $\ch^{-1}$ exists on $\{\psi_0\}^\perp$ and since $\psi_0'\in \{\psi_0\}^\perp$, it 
 is possible to compute $ \ch^{-1}\psi_0'$. 

This could be done almost explicitly, in any case explicit enough, so that we may compute $\dpr{\ch^{-1}\psi_0'}{\psi_0'}$ precisely in terms of $c$. We have 
$$
\left|
\begin{array}{l}
-\mu^2 f''+f -\f{\vp^2}{2 \mu^2} f+ \vp g'=\vp'' \\
-\mu^2 g''-(\vp f)' = -(\f{\vp^2}{2\mu^2})'. 
\end{array}
\right.
$$
Integrating once in the second equation yields 
$$
g'=\f{1}{\mu^2}(\f{\vp^2}{2\mu^2}-\vp f).
$$
Plugging this in the first equation yields 
\begin{equation}
\label{k:10}
-\mu^2 f''+f- \f{3 \vp^2}{2\mu^2} f +\f{\vp^3}{2\mu^4}=\vp''.
\end{equation}
On the other hand, recalling equation \eqref{eq:k1} in the form 
\begin{equation}
\label{k:15}
-\mu^2 \vp''+\vp-\f{\vp^3}{2\mu^2}=0,
\end{equation}
take a derivative in the parameter $\mu$ in \eqref{k:15}. We get 
$$
-\mu^2 (\p_\mu \vp)''+ \p_\mu \vp- \f{3 \vp^2}{2\mu^2} \p_\mu \vp+\f{\vp^3}{\mu^3}=2\mu \vp''.
$$
Dividing by $2\mu$ yields the relation 
\begin{equation}
\label{k:20}
-\mu^2 \left(\f{\p_\mu \vp}{2\mu}\right)''+ \left(\f{\p_\mu \vp}{2\mu}\right)- \f{3 \vp^2}{2\mu^2} \left(\f{\p_\mu \vp}{2\mu}\right) +\f{\vp^3}{\mu^3}=  \vp''.
\end{equation}
Comparing the equations \eqref{k:20} and \eqref{k:10} clearly implies (since we know that there is an unique solution by the invertibility of $\ch$ on $\{\psi_0\}^\perp$)  that 
$$
f=\f{\p_\mu \vp}{2\mu}.
$$
Thus, it follows that 
$$
g'=\f{1}{\mu^2}(\f{\vp^2}{2\mu^2}-\vp f)= g'=\f{1}{\mu^2}(\f{\vp^2}{2\mu^2}-\vp \f{\p_\mu \vp}{2\mu} )
$$
Now, 
\begin{eqnarray*}
& & \dpr{\ch \left(\begin{array}{c} \vp'' \\ -(\f{\vp^2}{2\mu^2})' \end{array}\right)}{\left(\begin{array}{c} \vp'' \\ -(\f{\vp^2}{2\mu^2})' \end{array}\right)}=
\dpr{\left(\begin{array}{c} f \\g \end{array}\right)}{\left(\begin{array}{c} \vp'' \\ -(\f{\vp^2}{2\mu^2})' \end{array}\right)}=  
  \dpr{\f{\p_\mu \vp}{2\mu}}{\vp''}+\dpr{g'}{\f{\vp^2}{2\mu^2}}= \\
  & & = -\f{\p_\mu[\int |\vp'(y)|^2dy]}{4\mu} -
  \f{\p_\mu[\int \vp^4(y) dy]}{16 \mu^5}+\f{\int \vp^4(y) dy}{4\mu^6}.
\end{eqnarray*}
Now since 
\begin{eqnarray*}
& & \int |\vp'(y)|^2dy = 4\int (sech'(y/\mu))^2 dy=4\mu \int (sech'(x))^2 dx\\
& & \int \vp^4(y) dy= 16\mu^5 \int sech(x)^4 dx,
\end{eqnarray*}
we have 
$$
\dpr{\ch \left(\begin{array}{c} \vp'' \\ -(\f{\vp^2}{2\mu^2})' \end{array}\right)}{\left(\begin{array}{c} \vp'' \\ -(\f{\vp^2}{2\mu^2})' \end{array}\right)}=
-\f{1}{\mu}\left(\int (sech'(x))^2 dx+\int sech(x)^4 dx\right). 
$$
Recall however, that we have to also compute the normalization factor $m^2$. We have 
\begin{eqnarray*}
 \f{1}{m^2} &=&  \left\|\left(\begin{array}{c} \vp' \\ -\f{\vp^2}{2\mu^2} \end{array}\right)\right\|^2=\int |\vp'(y)|^2 dy+ \f{1}{4\mu^4} \int \vp^4(y)dy= \\
& =& 
4\mu(\left(\int (sech'(x))^2 dx+\int sech(x)^4 dx\right).
\end{eqnarray*}
Thus, putting the last two formulas together 
\begin{equation}
\label{k:30}
\dpr{\ch^{-1}\psi_0'}{\psi_0'}=m^2 \dpr{\ch \left(\begin{array}{c} \vp'' \\ -(\f{\vp^2}{2\mu^2})' \end{array}\right)}{\left(\begin{array}{c} \vp'' \\ -(\f{\vp^2}{2\mu^2})' \end{array}\right)}=-\f{1}{4\mu^2}. 
\end{equation}
From this identity, we get that $\dpr{\ch^{-1}\psi_0'}{\psi_0'}<0$ for all values of the parameters, so in particular we claim that \eqref{B} is satisfied. 

Indeed, if we assume that $0=\dpr{\phi'}{\psi_0}=-\dpr{\phi}{\psi_0'}$, it would follow that $\psi_0'\in \{\phi\}^\perp$. Since $\ch|_{\{\phi\}^\perp}\geq 0$, it follows that 
$$
\dpr{\ch^{-1} \psi_0'}{\psi_0'}\geq 0,
$$
which is a contradiction with \eqref{k:30}. Properties \eqref{E} and \eqref{G} are obvious and hence, we may apply Theorem \ref{theo:5}. Since
$$
\om^*(\ch)=\f{1}{2\sqrt{-\dpr{\ch^{-1}\psi_0'}{\psi_0'}}}=\mu, 
$$
we have linear stability for all $c$ satisfying 
$
1>|c|\geq \mu=\sqrt{1-c^2},
$
and instability for $0\leq |c|< \mu=\sqrt{1-c^2}$. Solving the inequalities yields stability for $c: |c|\in [\f{\sqrt{2}}{2},1)$ and instability for $c: |c|\in 
[0,\f{\sqrt{2}}{2})$. Theorem \ref{theo:KGZ} is established. 

\section{Proof of Theorem \ref{theo:beam}} 
\label{sec:beam} 
Since we cannot verify theoretically the assumptions\footnote{but, as we have alluded to before, they hold nice and steady in numerical simulations, \cite{DS}} of Theorem \ref{theo:5}, we have put the appropriate assumptions, essentially requiring them, so that Theorem \ref{theo:5} applies.  

First, the case $c=0$  always presents instability, since $\ch_c$ has a negative eigenvalue. So, assume $c>0$, the cases $c<0$ being symmetric.  

Next, we need to compute $\ch^{-1}[\vp_c'']$. This is done by just taking a derivative with respect to the parameter $c$ in the defining equation \eqref{b:6}. We have 
$$
c^2(\p_c \vp)''+2c\vp_c''+(\p_c \vp)''''+(\p_c \vp)-p\vp^{p-1} (\p_c \vp) =0
$$
or equivalently $\ch[\p_c \vp]= -2c\vp_c''$. Assuming $c\neq 0$, we conclude 
$\ch^{-1}[\vp_c'']=-\f{1}{2c} \p_c \vp$. Now, 
$$
\dpr{\ch^{-1}\psi_0'}{\psi_0'}=-\f{1}{2c\|\vp'\|^2} \dpr{\p_c \vp_c}{\vp''}=
\f{1}{2c\|\vp'\|^2} \dpr{\p_c \vp'}{\vp'}=\f{\p_c(\|\vp_c'\|^2)}{4c\|\vp'\|^2}.
$$
According to Theorem \ref{theo:5}, we have instability for all $c$ such that 
$$
\p_c\|\vp'_c\|\geq 0
$$
and moreover, if $\p_c\|\vp_c'\|<0$, we have instability for all $c>0$, satisfying 
$$
c< \f{1}{2 \sqrt{-\f{\p_c \|\vp'_c\|^2}{4c\|\vp'_c\|^2} }}
$$
or 
$$
\sqrt{c}<\f{1}{\sqrt{-\f{\p_c \|\vp'_c\|^2}{\|\vp'_c\|^2}}}.
$$
In the complementary range, 
we have spectral stability.

\appendix

\section{Proof of Lemma \ref{cont}}   
  The term $\f{\la^2-\de^2}{4\om^2\la^2}$ is clearly continuous in both $\om, \la$, so we concentrate on the continuity of the mapping $(\om, \la)\to 
  \dpr{[H+\la^2+2\om\la P_0\p_x P_0]^{-1}[\phi']}{\phi'}$.

 We now need to show the continuity of the map stated above. 
 Taking a sequence $(\om_n, \la_n)\to (\om_0,\la_0)$ and denoting 
 \begin{eqnarray*}
 R_n &=& (H+\la_n^2+2\om_n\la_n P_0\p_x P_0)^{-1} \\
 R_0 &=& (H+\la_0^2+2\om_0\la_0 P_0\p_x P_0)^{-1},
 \end{eqnarray*}
   we need to show that 
$$
 \dpr{R_n[\phi']}{\phi'}\to \dpr{R_0[\phi']}{\phi'},
 $$
  which follows from   
  $$
  \|(R_n-R_0)P_0\|_{L^2\to L^2}\to 0,
  $$
 which remains to be proved. 
  By the resolvent identity, we have 
  $$
  R_n-R_0=-R_n(\la_n^2-\la_0^2+2(\om_n\la_n-\om_0\la_0)P_0\p_xP_0) R_0
  $$
  and since $|\la_n^2-\la_0^2+2(\om_n\la_n-\om_0\la_0)|
  \leq C(|\om_n|+|\om_0|+|\la_0|+|\la_n|)(|\om_n-\om_0|+|\la_n-\la_0|)$, it will 
  suffice to prove 
  \begin{equation}
  \label{m:110}
  \limsup_n \|R_n R_0 P_0\|_{\cb(L^2)}+\|R_n (P_0\p_xP_0) R_0 P_0\|_{\cb(L^2)}<\infty. 
  \end{equation}
  The first estimate follows from Proposition \ref{prop:5} 
  $$
  \limsup_n\|R_n R_0 P_0\|_{\cb(L^2)}\leq 
  \limsup_n\|R_n P_0\|_{\cb(L^2)} \|R_0 P_0\|_{\cb(L^2)}\leq 
  \limsup_n (\la_n^{-2} \la_0^{-2})=\la_0^{-4}.
  $$
  For the second term, we further use the resolvent identity to write   
  $$
  R_n (P_0\p_xP_0) R_0 P_0=  (H+\la_n^2)^{-1}  P_0 \p_x R_0 P_0- 
  2\om_n \la_n R_n P_0\p_x P_0 (H+\la_n^2)^{-1} \p_x R_0 P_0,
  $$
  and we do similar expansion on the right-side with $R_0=(H+\la_0^2)^{-1}+...$. 
  Clearly, in the above formula, we can pair the operators $\p_x$ with 
  various resolvents of the form  $(H+\la_n^2)^{-1}, (H+\la_0^2)^{-1}$, through \eqref{E} to obtain the desired boundedness results.

\end{document}